\documentclass[11pt,reqno]{amsart}

\usepackage{amssymb, amsmath, amsthm}
\usepackage[backref]{hyperref}
\usepackage[alphabetic,backrefs,lite]{amsrefs}
\usepackage{dlfltxbcodetips}
\usepackage{verbatim}
\usepackage{amscd}   
\usepackage[all]{xy} 
\usepackage{youngtab} 
\usepackage{young} 
\usepackage{tikz}

\renewcommand{\a}{\alpha}

\newcommand{\A}{\textit{\textbf{u}}}
\renewcommand{\b}{\beta}

\renewcommand{\c}{\gamma}
\renewcommand{\d}{\underline{d}}

\renewcommand{\k}{\kappa}
\renewcommand{\ll}{\lambda}

\renewcommand{\o}{\otimes}

\renewcommand{\r}{\underline{r}}
\newcommand{\s}{\sigma}

\renewcommand{\sl}{\mathfrak{sl}}

\renewcommand{\t}{\tau}

\newcommand{\GL}{{GL}}
\newcommand{\Hom}{\operatorname{Hom}} 
\newcommand{\Sym}{\operatorname{Sym}} 

\newcommand{\defi}[1]{{\upshape\sffamily #1}}

\newcommand{\bb}[1]{\mathbb{#1}}
\renewcommand{\rm}[1]{\mathrm{#1}}
\newcommand{\mc}[1]{\mathcal{#1}}

\newcommand{\ccircle}[1]{*+<1ex>[o][F-]{#1}}
\newcommand{\ccirc}[1]{\xymatrix@1{*+<1ex>[o][F-]{#1}}}

\def\ccone{\ccirc{1}}
\def\cctwo{\ccirc{2}}
\def\ccthree{\ccirc{3}}

\def\cck{\ccirc{k}}
\def\ccx{\ccirc{x}}
\def\ccy{\ccirc{y}}
\def\ccz{\ccirc{z}}

\def\kplone{k+1}

\def\aone{a_1}
\def\ai{a_i}
\def\at{a_t}
\def\ccaone{\ccirc{\aone}}
\def\ccai{\ccirc{\ai}}
\def\ccat{\ccirc{\at}}
\def\ccb{\ccirc{b}}

\def\lra{\longrightarrow}

\newtheorem{theorem}{Theorem}[section]
\newtheorem{lemma}[theorem]{Lemma}
\newtheorem{conjecture}[theorem]{Conjecture}
\newtheorem{proposition}[theorem]{Proposition}
\newtheorem{corollary}[theorem]{Corollary}
\newtheorem*{problem*}{Problem}

\newtheorem*{tmain}{Theorem \ref{thm:main}}

\theoremstyle{definition}
\newtheorem{definition}[theorem]{Definition}
\newtheorem*{definition*}{Definition}
\newtheorem{example}[theorem]{Example}

\theoremstyle{remark}
\newtheorem{remark}[theorem]{Remark}
\newtheorem*{remark*}{Remark}

\numberwithin{equation}{section}

\begin{document}

\title{Secant Varieties of Segre--Veronese Varieties}

\author{Claudiu Raicu}
\address{Department of Mathematics, Princeton University, Princeton, NJ 08544-1000\newline
\indent Institute of Mathematics ``Simion Stoilow" of the Romanian Academy}
\email{craicu@math.princeton.edu}


\date{\today}


\begin{abstract} We prove that the ideal of the variety of secant lines to a Segre--Veronese variety is generated in degree three by minors of flattenings. In the special case of a Segre variety this was conjectured by Garcia, Stillman and Sturmfels, inspired by work on algebraic statistics, as well as by Pachter and Sturmfels, inspired by work on phylogenetic inference. In addition, we describe the decomposition of the coordinate ring of the secant line variety of a Segre--Veronese variety into a sum of irreducible representations under the natural action of a product of general linear groups.
\end{abstract}

\maketitle

\section{Introduction}

The spaces of matrices (or $2$--tensors) are stratified according to rank by the secant varieties of Segre products of two projective spaces. The defining ideals of these secant varieties are known to be generated by minors of generic matrices. It is an important problem, with applications in algebraic statistics, biology, signal processing, complexity theory etc., to understand (border) rank varieties of higher order tensors. These are (upon taking closure) the classical secant varieties to Segre varieties, whose equations are far from being understood. To get an idea about the boundary of our knowledge, note that the Salmon problem \cite{allman}, which asks for the generators of the ideal of $\s_4(\bb{P}^3\times\bb{P}^3\times\bb{P}^3)$, the variety of secant $3$--planes to the Segre product of three projective $3$--spaces, is still unsolved (although its set--theoretic version has been recently resolved in \cites{friedland,friedland-gross}; see also \cite{bates-oeding}).

Flattenings (see Section \ref{sec:flattenings}) provide an easy tool for obtaining some equations for secant varieties of Segre products, but they are not sufficient in general, as can be seen for example in the case of the Salmon problem. Inspired by the study of Bayesian networks, Garcia, Stillman and Sturmfels conjectured \cite[Conjecture 21]{GSS} that flattenings give all the equations of the first secant variety of the Segre variety. This conjecture also appeared at the same time in a biological context, namely in work of Pachter and Sturmfels on phylogenetic inference \cite[Conjecture 13]{pachter-sturmfels}.

\begin{conjecture}[Garcia--Stillman--Sturmfels]\label{conj:gss} The ideal of the secant line variety of a Segre product of projective spaces is generated by $3\times 3$ minors of flattenings.
\end{conjecture}

The set--theoretic version of this conjecture was obtained by Landsberg and Manivel \cite{lan-manivel}, as well as the case of a $3$--factor Segre product. The $2$--factor case is classical, while the $4$--factor case was resolved by Landsberg and Weyman \cite{lan-weyman}. The $5$--factor case was proved by Allman and Rhodes \cite{allman-rhodes}. We prove the GSS conjecture in Corollary \ref{cor:gss} as a consequence of our main result, Theorem \ref{thm:main}, which is the corresponding statement for Segre--Veronese varieties.

It is a general fact that for a subvariety $X$ in projective space which is not contained in a hyperplane, the ideal of the variety $\s_k(X)$ of secant $(k-1)$--planes to $X$ has no equations in degree less than $k+1$. If $X=G/P$ is a rational homogeneous variety, a theorem of Kostant (see \cite{landsberg}) states that the ideal of $X$ is generated in the smallest possible degree (i.e. in degree two), and Landsberg and Manivel asked whether this is also true for the first secant variety of $X$ \cite{lan-manivel}. It turns out that when $X$ is the $D_7$--spinor variety, there are in fact no cubics in the ideal of $\s_2(X)$ (see \cite{lan-weyman-chss} and \cite{manivel}). In Theorem \ref{thm:main}, we provide a class of $G/P$'s, the Segre--Veronese varieties for which the answer to the question of Landsberg and Manivel is positive. We obtain furthermore an explicit decomposition into irreducible representations of the homogeneous coordinate ring of the secant line variety of a Segre--Veronese variety, thus making it possible to compute the Hilbert function for this class of varieties. This can be regarded as a generalization of the computation of the degree of these secant varieties in \cite{cox-sidman}.

Before stating the main theorem, we establish some notation. For a vector space $V$, $V^*$ denotes its dual, and $\bb{P}V$ denotes the projective space of lines in $V$. If $\mu=(\mu_1\geq\mu_2\geq\cdots)$ is a partition, $S_{\mu}$ denotes the corresponding Schur functor (if $\mu_2=0$ we get symmetric powers, whereas if all $\mu_i=1$, we get exterior powers). For positive integers $d_1,\cdots,d_n$, $SV_{d_1,\cdots,d_n}$ denotes the Segre--Veronese embedding of a product of $n$ projective spaces via the complete linear system of the ample line bundle $\mc{O}(d_1,\cdots,d_n)$. $\s_2(X)$ denotes the variety of secant lines to $X$.

\begin{tmain}
 Let $X=SV_{d_1,\cdots,d_n}(\bb{P}V_1^*\times\bb{P}V_2^*\times\cdots\times\bb{P}V_n^*)$ be a Segre--Veronese variety, where each $V_i$ is a vector space of dimension at least $2$ over a field $K$ of characteristic zero. The ideal of $\s_2(X)$ is generated by $3\times 3$ minors of flattenings, and moreover, for every nonnegative integer $r$ we have the decomposition of the degree $r$ part of its homogeneous coordinate ring
\[K[\s_2(X)]_r=\bigoplus_{\substack{\ll=(\ll^1,\cdots,\ll^n)\\ \ll^i\vdash rd_i}}(S_{\ll^1}V_1\o\cdots\o S_{\ll^n}V_n)^{m_{\ll}},\]
where $m_{\ll}$ is obtained as follows. Set 
\[f_{\ll}=\max_{i=1,\cdots,n}\left\lceil\frac{\ll_2^i}{d_i}\right\rceil,\quad e_{\ll}=\ll^1_2+\cdots+\ll^n_2.\]
If some partition $\ll^i$ has more than two parts, or if $e_{\ll}<2f_{\ll}$, then $m_{\ll}=0$. If $e_{\ll}\geq r-1$, then $m_{\ll}=\lfloor r/2\rfloor-f_{\ll}+1$, unless $e_{\ll}$ is odd and $r$ is even, in which case $m_{\ll}=\lfloor r/2\rfloor-f_{\ll}$. If $e_{\ll}<r-1$ and $e_{\ll}\geq 2f_{\ll}$, then $m_{\ll}=\lfloor (e_{\ll}+1)/2\rfloor-f_{\ll}+1$, unless $e_{\ll}$ is odd, in which case  $m_{\ll}=\lfloor (e_{\ll}+1)/2\rfloor-f_{\ll}$.
\end{tmain}

Theorem \ref{thm:main} has further consequences to deriving certain plethystic formulas for decomposing (in special cases) symmetric powers of triple tensor products (Corollary \ref{cor:pletysm3factors}a)) and Schur functors applied to tensor products of two vector spaces (Corollary \ref{cor:pletysm3factors}b)), or even symmetric pletyhsm (Corollary \ref{cor:symmetricplethysm}).

The main technique introduced in this work does not seek to employ the particularities of specific instances of Segre--Veronese varieties, but instead tries to capture only the essential features that are shared between all these varieties. We work in some sense with spaces of ``generic tensors'', and rather concentrate on their ``generic equations''. The latter are representations of products of symmetric groups, which can be defined abstractly with no relation to spaces of tensors (although what led us to them was their realization as zero--weight spaces of particular tensor representations). The main point is that the generic equations yield, by a process of \defi{specialization}, the equations of any specific secant variety of a Segre--Veronese variety. One can also go back, via \defi{polarization}, from the equations of a specific secant variety to (a subset of) the generic equations. The main tools that we employ in analyzing the generic equations of the varieties of secant lines are combinatorial: graph theory and tableaux combinatorics. We hope that similar methods, particularly replacing graph theory with the theory of simplicial complexes, could be used to give an analogous picture for higher secant varieties. The main goal of our work is to set up a general framework that would help understand arbitrary secant varieties, and illustrate how insights from combinatorics occur naturally in this framework, providing new results in the case of the varieties of secant lines.

Finding equations for higher secant varieties of Segre--Veronese varieties turns out to be a delicate task, even in the case of two factors $(n=2)$ with not too positive embeddings (small $d_1,d_2$). Recent progress in this direction has been obtained by Cartwright, Erman and Oeding \cite{ceo}.

Since finding precise descriptions of the equations, and more generally syzygies, of secant varieties to Segre--Veronese varieties constitutes such an intricate project, much of the current effort is directed to finding more qualitative statements. Draisma and Kuttler \cite{draisma-kuttler} prove that for each $k$, there is an uniform bound $d(k)$ such that the $(k-1)$--st secant variety of any Segre variety is cut out (set--theoretically) by equations of degree at most $d(k)$. Theorem \ref{thm:main} implies that $d(2)=3$, even ideal theoretically.

For higher syzygies, Snowden proves that all the syzygies of Segre varieties are obtained from a finite amount of data via an iterative process \cite{snowden}. It would be interesting to know if the same result holds for the secant varieties. This would generalize the result of Draisma and Kuttler. For Veronese varieties, the asymptotic picture of the Betti tables is described in work of Ein and Lazarsfeld \cite{ein-lazarsfeld}. Again, it would be desirable to have analogous results for secant varieties.

The structure of the paper is as follows. In Section~\ref{chap:prelim} we give the basic definitions for secant varieties and Segre--Veronese varieties. We introduce the basic notions from Representation Theory that are used throughout the work, and describe the process of flattening a tensor, which leads to the notion of a flattening matrix. Section~\ref{chap:equations} builds the framework for analyzing the equations and homogeneous coordinate rings of arbitrary secant varieties of Segre--Veronese varieties. Even though we were only able to work out the details of this analysis in the case of the first secant variety, we believe that the general method of approach may be used to shed some light on the case of higher secant varieties. In particular, the new insight of concentrating on the ``generic equations'' is presented in detail and in the generality needed to deal with arbitrary secant varieties. Section~\ref{chap:secantline} is inspired by a conjecture of Garcia, Stillman and Sturmfels, describing the generators of the ideal of the variety of secant lines to a Segre variety. We prove more generally that this description holds for the first secant variety of a Segre--Veronese variety. We also give a representation theoretic decomposition of the coordinate ring of this variety, which allows us to deduce certain plethystic formulas based on known computations of dimensions of secant varieties of Segre varieties.

\section{Preliminaries}\label{chap:prelim}

Throughout this work, $K$ denotes a field of characteristic $0$. All the varieties we study are of finite type over $K$, and are reduced and irreducible. $\bb{P}^N$ denotes the $N$--dimensional projective space over $K$. We write $\bb{P}W$ for $\bb{P}^N$ when we think of $\bb{P}^N$ as the space of $1$--dimensional subspaces (lines) in a vector space $W$ of dimension $N+1$ over $K$. Given a nonzero vector $w\in W$, we denote by $[w]$ the corresponding line. The coordinate ring of $\bb{P}W$ is $\rm{Sym}(W^*)$, the symmetric algebra on the vector space $W^*$ of linear functionals on $W$.

\subsection{Secant Varieties}\label{sec:secants}

\begin{definition} Given a subvariety $X\subset\bb{P}^N$, the \defi{$(k-1)$-st secant variety} of $X$, denoted $\s_k(X)$, is the closure of the union of linear subspaces spanned by $k$ points on $X$:
\[\s_k(X)=\overline{\bigcup_{x_1,\cdots,x_k\in X}\bb{P}_{x_1,\cdots,x_k}}.\]
\end{definition}

Alternatively, if we write $\bb{P}^N=\bb{P}W$ for some vector space $W$, and let $\hat{X}\subset W$ denote the cone over $X$, then we can define $\s_k(X)$ by specifying its cone $\widehat{\s_k(X)}$. This is the closure of the image of the map
\[s:\hat{X}\times\cdots\times\hat{X}\lra W,\]
defined by
\[s(x_1,\cdots,x_k)=x_1+\cdots+x_k.\]

The main problem we are concerned with is

\begin{problem*} Given (the equations of) $X$, determine (the equations of) $\s_k(X)$.
\end{problem*}

More precisely, given the homogeneous ideal $I(X)$ of the subvariety $X\subset\bb{P}W$, we would like to describe the generators of $I(\s_k(X))$. Alternatively, we would like to understand the homogeneous coordinate ring of $\s_k(X)$, which we denote by $K[\s_k(X)]$. As we will see, this is a difficult problem even in the case when $X$ is simple, i.e. isomorphic to a projective space (or a product of such). There is thus little hope of giving an uniform satisfactory answer in the generality we posed the problem. However, the following observation provides a general approach to the problem, which we exploit in the future sections.

The ideal/homogeneous coordinate ring of a subvariety $Y\subset\bb{P}W$ coincides with the ideal/affine coordinate ring of its cone $\hat{Y}\subset W$, hence our problem is equivalent to understanding $I(\widehat{\s_k(X)})$ and $K[\widehat{\s_k(X)}]$. The morphism $s$ of affine varieties defined above corresponds to a ring map
\[s^{\#}:\rm{Sym}(W^*)\to K[\hat{X}\times\cdots\times\hat{X}]=K[\hat{X}]\o\cdots\o K[\hat{X}].\]
We have that $I(\widehat{\s_k(X)})$ and $K[\widehat{\s_k(X)}]$ are the kernel and image respectively of $s^{\#}$. The main focus for us will be on the case when $X$ is a Segre-Veronese variety (described in the following section), and $k=2$.

\subsection{Segre-Veronese Varieties}\label{sec:segreveronese}

Consider vector spaces $V_1,\cdots,V_n$ of dimensions $m_1,\cdots,m_n\geq 2$ respectively, with duals $V_1^*,\cdots,V_n^*$, and positive integers $d_1,\cdots,d_n$. We let 
\[X=\bb{P}V_1^*\times\cdots\times\bb{P}V_n^*\]
and think of it as a subvariety in projective space via the embedding determined by the line bundle $\mc{O}_X(d_1,\cdots,d_n)$. Explicitly, $X$ is the image of the map
\[SV_{d_1,\cdots,d_n}:\bb{P}V_1^*\times\cdots\times\bb{P}V_n^*\to\bb{P}(\Sym^{d_1}V_1^*\o\cdots\o\Sym^{d_n}V_n^*)\]
given by
\[([e_1],\cdots,[e_n])\mapsto[e_1^{d_1}\o\cdots\o e_n^{d_n}].\]
We call $X$ a \defi{Segre-Veronese variety}.

For such $X$ we prove that $I(\s_2(X))$ is generated in degree $3$ and we describe the decomposition of $K[\s_2(X)]$ into a sum of irreducible representations of the product of general linear groups $\GL(V_1)\times\cdots\times\GL(V_n)$ (Theorem \ref{thm:main}).

When $n=1$ we set $d=d_1$, $V=V_1$. The image of $SV_d$ is the $d$-th Veronese embedding, or $d$-uple embedding of the projective space $\bb{P}V^*$, which we denote by $\rm{Ver}_d(\bb{P}V^*)$. When $d_1=\cdots=d_n=1$, the image of $SV_{1,1,\cdots,1}$ is the Segre variety $\rm{Seg}(\bb{P}V_1^*\times\cdots\times\bb{P}V_n^*)$. An element of $\Sym^{d_1}V_1^*\o\cdots\o\Sym^{d_n}V_n^*$ is called a \defi{(partially symmetric) tensor}. The points in the cone over the Segre-Veronese variety are called \defi{pure tensors}. 

\subsection{Representation Theory}\label{sec:reptheory}

We refer the reader to \cite{ful-har} for the basic representation theory of symmetric and general linear groups. Given a positive integer $r$, a partition $\mu$ of $r$ is a nonincreasing sequence of nonnegative integers $\mu_1\geq\mu_2\geq\cdots$ with $r=\sum\mu_i$. We write $\mu=(\mu_1,\mu_2,\cdots)$. Alternatively, if $\mu$ is a partition having $i_j$ parts equal to $\mu_j$ for all $j$, then we write $\mu=(\mu_1^{i_1}\mu_2^{i_2}\cdots)$. To a partition $\mu=(\mu_1,\mu_2,\cdots)$ we associate a \defi{Young diagram} which consists of left-justified rows of boxes, with $\mu_i$ boxes in the $i$-th row. For $\mu=(5,2,1)$, the corresponding Young diagram is
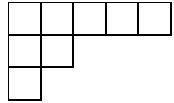
\[\yng(5,2,1)\]

For a vector space $W$, a positive integer $r$ and a partition $\mu$ of $r$, we denote by $S_{\mu}W$ the corresponding irreducible representation of $\GL(W)$: $S_{\mu}$ are commonly known as \defi{Schur functors}, and we make the convention that $S_{(d)}$ denotes the symmetric power functor, while $S_{(1^d)}$ denotes the exterior power functor. We write $S_r$ for the symmetric group on $r$ letters, and $[\mu]$ for the irreducible $S_r$-representation corresponding to $\mu$: $[(d)]$ denotes the trivial representation and $[(1^d)]$ denotes the sign representation.

Given a positive integer $n$ and a sequence of nonnegative integers $\r=(r_1,\cdots,r_n)$, we define an \defi{$n$-partition} of $\r$ to be an $n$-tuple of partitions $\ll=(\ll^1,\cdots,\ll^n)$, with $\ll^j$ partition of $r_j$, $j=1,\cdots,n$. We write $\ll^j\vdash r_j$ and $\ll\vdash^n\r$. Given vector spaces $V_1,\cdots,V_n$ as above, we often write $\GL(V)$ for $\GL(V_1)\times\cdots\times\GL(V_n)$. We write $S_{\ll}V$ for the irreducible $\GL(V)$-representation $S_{\ll^1}V_1\o\cdots\o S_{\ll^n}V_n$. Similarly, we write $[\ll]$ for the irreducible representation $[\ll^1]\o\cdots\o [\ll^n]$ of the $n$-fold product of symmetric groups $S_{\r}=S_{r_1}\times\cdots\times S_{r_n}$. We have

\begin{lemma}[Schur-Weyl duality]\label{lem:SchurWeyl}
\[V_1^{\o r_1}\o\cdots\o V_n^{\o r_n}=\bigoplus_{\ll\vdash^n\r}[\ll]\o S_{\ll}V.\]
\end{lemma}

Most of the group actions we consider are left actions, denoted by $\cdot$. We use the symbol $*$ for right actions, to distinguish them from left actions.

For a subgroup $H\subset G$ and representations $U$ of $H$ and $W$ of $G$, we write
\[\rm{Ind}_H^G(U)=K[G]\o_{K[H]} U,\ \rm{ and }\ \rm{Res}_H^G(W)=W_H,\]
for the \defi{induced representation} of $U$ and \defi{restricted representation} of $W$, where $K[M]$ denotes the group algebra of a group $M$, and $W_H$ is just $W$, regarded as an $H$-module. We write $W^G$ for the \defi{$G$-invariants} of the representation $W$, i.e.
\[W^G=\rm{Hom}_G({\bf 1},W)\subset \rm{Hom}_K({\bf 1},W)=W,\]
where $\bf 1$ denotes the trivial representation of $G$. 

\begin{remark}\label{rem:invariants} If $G$ is finite, let
\[s_G=\sum_{g\in G}g\in K[G].\]
We can realize $W^G$ as the image of the map $W\lra W$ given by
\[
 w\mapsto s_G\cdot w.
\]
Assume furthermore that $H\subset G$ is a subgroup, and let $s_H$ denote the corresponding element in $K[H]$. We have a natural inclusion of the trivial representation of $H$
\[{\bf 1}\hookrightarrow K[H],\quad 1\mapsto s_H,\]
which after tensoring with $K[G]$ becomes
\[\rm{Ind}_H^G({\bf 1})=K[G]\o_{K[H]}{\bf 1}\hookrightarrow K[G]\o_{K[H]}K[H]\simeq K[G],\]
so that we can identify $\rm{Ind}_H^G({\bf 1})$ with $K[G]\cdot s_H$.
\end{remark}
We have

\begin{lemma}[Frobenius reciprocity]\label{lem:Frobenius}
\[W^H=\rm{Hom}_H({\bf 1},\rm{Res}_H^G(W))=\rm{Hom}_G(\rm{Ind}_H^G({\bf 1}),W).\]
\end{lemma}

Given an $n$-partition $\ll=(\ll^1,\cdots,\ll^n)$ of $\r$, we define an \defi{$n$-tableau} of shape $\ll$ to be an $n$-tuple $T=(T^1,\cdots,T^n)$, which we usually write as $T^1\o\cdots\o T^n$, where each $T^i$ is a tableau of shape $\ll^i$. A tableau is \defi{canonical} if its entries index its boxes consecutively from left to right, and top to bottom. We say that $T$ is canonical if each $T^i$ is, in which case we write $T_{\ll}$ for $T$. If $T=(\ll^1,\ll^2)$, with $\ll^1=(3,2)$, $\ll^2=(3,1,1)$, then the canonical $2$-tableau of shape $\ll$ is
\[\Yvcentermath1 \young(123,45)\o\young(123,4,5).\]

We consider the subgroups of $S_{\r}$
\[R_{\ll} = \{g\in S_{\r}:g\ \textrm{ preserves each row of }\ T_{\ll}\},\]
\[C_{\ll} = \{g\in S_{\r}:g\ \textrm{ preserves each column of }\ T_{\ll}\}\]
and define the symmetrizers
\[a_{\ll}=\sum_{g\in R_{\ll}}g,\quad b_{\ll}=\sum_{g\in C_{\ll}}\rm{sgn}(g)\cdot g,\quad c_{\ll}=a_{\ll}\cdot b_{\ll},\]
with $\rm{sgn}(g)=\prod_i\rm{sgn}(g_i)$ for $g=(g_1,\cdots,g_n)\in S_{\r}$, where $\rm{sgn}(g_i)$ denotes the signature of the permutation $g_i$.

The $\GL(V)$- (or $S_{\r}$-) representations $W$ that we consider decompose as a direct sum of $S_{\ll}V$'s (or $[\ll]$'s) with $\ll\vdash^n\r$. We write
\[W=\bigoplus_{\ll}W_{\ll},\]
where $W_{\ll}\simeq (S_{\ll}V)^{m_{\ll}}$ (or $W_{\ll}\simeq [\ll]^{m_{\ll}}$) for some nonnegative integer $m_{\ll}=m_{\ll}(W)$, called the \defi{multiplicity} of $S_{\ll}V$ (or $[\ll]$) in $W$. We call $W_{\ll}$ the \defi{$\ll$-part} of the representation $W$.

Recall that $m_j$ denotes the dimension of $V_j$, $j=1,\cdots,n$. We fix bases 
\[\mc{B}_j=\{x_{ij}:i=1,\cdots,m_j\}\]
for $V_j$ ordered by $x_{ij}>x_{i+1,j}$. We choose the \defi{maximal torus} $T=T_1\times\cdots\times T_n\subset\GL(V)$, with $T_j$ being the set of diagonal matrices with respect to $\mc{B}_j$. We choose the \defi{Borel subgroup} of $\GL(V)$ to be $B=B_1\times\cdots\times B_n$, where $B_j$ is the subgroup of upper triangular matrices in $\GL(V_j)$ with respect to $\mc{B}_j$. Given a $\GL(V)$-representation $W$, a \defi{weight vector} $w$ with \defi{weight} $a=(a_1,\cdots,a_n)$, $a_i\in T_i^*$, is a nonzero vector in $W$ with the property that for any $t=(t_1,\cdots,t_n)\in T$,
\[t\cdot w=a_1(t_1)\cdots a_n(t_n) w.\]
The vectors with this property form a vector space called the \defi{$a$-weight space} of $W$, which we denote by $\rm{wt}_a(W)$.

A \defi{highest weight vector} of a $\GL(V)$-representation $W$ is an element $w\in W$ invariant under $B$. $W=S_{\ll}V$ has a unique (up to scaling) highest weight vector $w$ with corresponding weight $\ll=(\ll^1,\cdots,\ll^n)$. In general, we define the \defi{$\ll$-highest weight space} of a $\GL(V)$-representation $W$ to be the set of highest weight vectors in $W$ with weight $\ll$, and denote it by $\rm{hwt}_{\ll}(W)$. If $W$ is an $S_{\r}$-representation, the \defi{$\ll$-highest weight space} of $W$ is the vector space $\rm{hwt}_{\ll}(W)=c_{\ll}\cdot W\subset W$, where $c_{\ll}$ is the Young symmetrizer defined above. In both cases, $\rm{hwt}_{\ll}(W)$ is a vector space of dimension $m_{\ll}(W)$.

\subsection{Flattenings}\label{sec:flattenings} 

Given decompositions $d_i=a_i+b_i$, with $a_i,b_i\geq 0$, $i=1,\cdots,n$, we let $A=(a_1,\cdots,a_n)$, $B=(b_1,\cdots,b_n)$, so that $\d=(d_1,\cdots,d_n)=A+B$, and embed
\[\Sym^{d_1}V_1^*\o\cdots\o\Sym^{d_n}V_n^*\hookrightarrow V_A^*\o V_B^*\]
in the usual way, where
\[V_A=\Sym^{a_1}V_1\o\cdots\o\Sym^{a_n}V_n,\quad V_B=\Sym^{b_1}V_1\o\cdots\o\Sym^{b_n}V_n.\]
This embedding allows us to \defi{flatten} any tensor in $\Sym^{d_1}V_1^*\o\cdots\o\Sym^{d_n}V_n^*$ to a $2$-tensor, i.e. a matrix, in $V_A^*\o V_B^*$. We call such a matrix an \defi{$(A,B)$-flattening} of our tensor. If $|A|=a_1+...+a_n$ then we also say that this matrix is an $|A|$-flattening, or a $|B|$-flattening, by symmetry.

We obtain an inclusion
\[SV_{d_1,\cdots,d_n}(\bb{P}V_1^*\times\cdots\times\bb{P}V_n^*)\hookrightarrow\rm{Seg}(\bb{P}V_A^*\times\bb{P}V_B^*),\]
and consequently
\[\s_k(SV_{d_1,\cdots,d_n}(\bb{P}V_1^*\times\cdots\times\bb{P}V_n^*))\hookrightarrow\s_k(\rm{Seg}(\bb{P}V_A^*\times\bb{P}V_B^*)),\]
where the latter secant variety coincides with (the projectivization of) the set of matrices of rank at most $k$ in $V_A^*\o V_B^*$. This set is cut out by the $(k+1)\times(k+1)$ minors of the generic matrix in $V_A^*\o V_B^*$. This observation yields equations for the secant varieties of Segre-Veronese varieties (see also \cite{landsberg}).

\begin{lemma}\label{lem:eqnflat} For any decomposition $\d=A+B$ and any $k\geq 1$, the ideal of $(k+1)\times(k+1)$ minors of the generic matrix given by the $(A,B)$-flattening of $\Sym^{d_1}V_1^*\o\cdots\o\Sym^{d_n}V_n^*$ is contained in the ideal of $\s_{k}(SV_{d_1,\cdots,d_n}(\bb{P}V_1^*\times\cdots\times\bb{P}V_n^*))$.
\end{lemma}

\begin{definition}\label{def:flattenings}
 We write $F_{A,B}^{k+1,r}(V)=F_{A,B}^{k+1,r}(V_1,\cdots,V_n)$ for the degree $r$ part of the ideal of $(k+1)\times(k+1)$ minors of the $(A,B)$-flattening.
\end{definition}

Note that the invariant way of writing the generators of the ideal of $(k+1)\times(k+1)$ minors of the $(A,B)$-flattening in the preceding lemma ($F_{A,B}^{k+1,k+1}(V)$) is as the image of the composition
\[\bigwedge^{k+1}V_A\o\bigwedge^{k+1}V_B\hookrightarrow\Sym^{k+1}(V_A\o V_B)\lra\Sym^{k+1}(\Sym^{d_1}V_1\o\cdots\o\Sym^{d_n}V_n),\]
where the first map is the usual inclusion map, while the last one is induced by the multiplication maps
\[\Sym^{a_i}V_i\o\Sym^{b_i}V_i\lra\Sym^{d_i}V_i.\]

\subsection{The ideal and coordinate ring of a Segre-Veronese variety}\label{sec:coordringSV}

If $X=SV_{d_1,\cdots,d_n}(\bb{P}V_1^*\times\cdots\times\bb{P}V_n^*)$, then the ideal $I(X)$ is generated by $2\times 2$ minors of flattenings, i.e. when $k=1$ the equations described in Lemma \ref{lem:eqnflat} are sufficient to generate the ideal of the corresponding variety. As for the homogeneous coordinate ring of a Segre-Veronese variety, we have the decomposition
\[K[X]=\bigoplus_{r\geq 0}(\Sym^{rd_1}V_1\o\cdots\o\Sym^{rd_n}V_n).\tag{*}\]
This decomposition will turn out to be useful in the next section, in conjunction with the map $s^{\#}$ defined in Section \ref{sec:secants}. In Section \ref{chap:secantline} we give a description of $K[\s_2(X)]$ analogous to (*), and prove that the $3\times 3$ minors of flattenings generate the homogeneous ideal of $\s_2(X)$.

The statements above regarding the ideal and coordinate ring of a Segre-Veronese variety hold more generally for rational homogeneous varieties $(G/P)$, and have been obtained in unpublished work by Kostant (see \cite{landsberg}).

\section{Equations of the secant varieties of a Segre-Veronese variety}\label{chap:equations}

This section introduces the main new tool for understanding the equations and coordinate rings of the secant varieties of Segre-Veronese varieties, from a representation theoretic/combinatorial perspective. All the subsequent work is based on the ideas described here. The usual method for analyzing the secant varieties of Segre-Veronese varieties is based on the representation theory of general linear groups. We review some of its basic ideas, including the notion of \defi{inheritance}, in Section \ref{sec:prolongations}. The new insight of restricting the analysis to special equations of the secant varieties, the ``generic equations'', is presented in Section \ref{sec:generic}. More precisely, we use Schur-Weyl duality to translate questions about representations of general linear groups into questions about representations of symmetric groups and tableaux combinatorics. The relationship between the two situations is made precise in Section \ref{sec:polarization}. One should think of the ``generic equations'' as a set of equations that give rise by specialization to all the equations of the secant varieties of Segre-Veronese varieties. Similarly, we have the ``generic flattenings'' which by specialization yield the usual flattenings.

\subsection{Multi-prolongations and inheritance}\label{sec:prolongations}

In this section $V_1,\cdots,V_n$ are (as always) vector spaces over a field $K$ of characteristic zero. We switch from the $\Sym^d$ notation to the more compact Schur functor notation $S_{(d)}$ described in Section \ref{sec:reptheory}. The homogeneous coordinate ring of $\bb{P}(S_{(d_1)}V_1^*\o\cdots\o S_{(d_n)}V_n^*)$ is 
\[S=\rm{Sym}(S_{(d_1)}V_1\o\cdots\o S_{(d_n)}V_n),\]
the symmetric algebra of the vector space $S_{(d_1)}V_1\o\cdots\o S_{(d_n)}V_n$. This vector space has a natural basis $\mc{B}=\mc{B}_{d_1,\cdots,d_n}$ consisting of tensor products of monomials in the elements of the bases $\mc{B}_1,\cdots,\mc{B}_n$ of $V_1,\cdots,V_n$. We write this basis, suggestively, as
\[\mc{B} = \Sym^{d_1}\mc{B}_1\o\cdots\o\Sym^{d_n}\mc{B}_n.\]
We can index the elements of $\mc{B}$ by $n$-tuples $\a=(\a_1,\cdots,\a_n)$ of multisets $\a_i$ of size $d_i$ with entries in $\{1,\cdots,m_i=\rm{dim}(V_i)\}$, as follows. The $\a$-th element of the basis $\mc{B}$ is
\[z_{\a} = (\prod_{i_1\in\a_1}x_{i_1,1})\o\cdots\o(\prod_{i_n\in\a_n}x_{i_n,n}),\]
and we think of $z_{\a}$ as a linear form in $S$.

We can therefore identify $S$ with the polynomial ring $K[z_{\a}]$. One would like to have a precise description of the ideal $I\subset S$ of polynomials vanishing on $\s_{k}(SV_{d_1,\cdots,d_n}(\bb{P}V_1^*\times\cdots\times\bb{P}V_n^*))$, but this is a very difficult problem, as mentioned in the introduction. We obtain such a description for the case $k=2$ in Theorem \ref{thm:main}. The case $k=1$ was already known, as described in Section \ref{sec:coordringSV}.

Given a positive integer $r$ and a partition $\mu=(\mu_1,\cdots,\mu_t)\vdash r$, we consider the set $\mc{P}_{\mu}$ of all (unordered) partitions of $\{1,\cdots,r\}$ of shape $\mu$, i.e.
\[\mc{P}_{\mu}=\left\{A=\{A_1,\cdots,A_t\}:|A_i|=\mu_i\ \rm{ and }\ \bigsqcup_{i=1}^t A_i=\{1,\cdots,r\}\right\},\]
as opposed to the set of ordered partitions where we take instead $A=(A_1,\cdots,A_t)$.


\begin{definition}\label{def:pimu} For a partition $\mu=(\mu_1^{i_1}\cdots\mu_s^{i_s})$ of $r$, we consider the map
\[\pi_{\mu}:S_{(r)}(S_{(d_1)}V_1\o\cdots\o S_{(d_n)}V_n)\lra \bigotimes_{j=1}^s S_{(i_j)}(S_{(\mu_j d_1)}V_1\o\cdots\o S_{(\mu_j d_n)}V_n),\]
given by
\[z_1\cdots z_r\mapsto\sum_{A\in\mc{P}_{\mu}}\bigotimes_{j=1}^s\prod_{\substack{B\in A\\|B|=\mu_j}} m(z_i:i\in B),\]
where 
\[m:(S_{(d_1)}V_1\o\cdots\o S_{(d_n)}V_n)^{\o\mu_j}\lra S_{(\mu_j d_1)}V_1\o\cdots\o S_{(\mu_j d_n)}V_n\]
denotes the usual componentwise multiplication map.

We write $\pi_{\mu}(V)$ or $\pi_{\mu}(V_1,\cdots,V_n)$ for the map $\pi_{\mu}$ just defined, when we want to distinguish it from its generic version (Definition \ref{def:pimugeneric}). We also write $U_r^{\d}(V)=U_r^{\d}(V_1,\cdots,V_n)$ and $U_{\mu}^{\d}(V)=U_{\mu}^{\d}(V_1,\cdots,V_n)$ for the source and target of $\pi_{\mu}(V)$ respectively (see Definitions \ref{def:urd} and \ref{def:upid} for the generic versions of these spaces).
\end{definition}

A more invariant way of stating Definition \ref{def:pimu} is as follows. If $\mu=(\mu_1,\cdots,\mu_t)$, then the map $\pi_{\mu}$ is the composition between the usual inclusion
\[S_{(r)}(S_{(d_1)}V_1\o\cdots\o S_{(d_n)}V_n)\hookrightarrow(S_{(d_1)}V_1\o\cdots\o S_{(d_n)}V_n)^{\o r}=\]
\[(S_{(d_1)}V_1\o\cdots\o S_{(d_n)}V_n)^{\o \mu_1}\o\cdots\o(S_{(d_1)}V_1\o\cdots\o S_{(d_n)}V_n)^{\o \mu_t}.\]
and the tensor product of the natural multiplication maps
\[m:(S_{(d_i)}V_i)^{\o\mu_j}\lra S_{(\mu_j d_i)}V_i.\]

\begin{example}\label{ex:pimu} Let $n=2$, $d_1=2$, $d_2=1$, $r=4$, $\mu=(2,2)=(2^2)$, $\rm{dim}(V_1)=2$, $\rm{dim}(V_2)=3$. Take
\[z_1=z_{(\{1,2\},\{1\})},\quad z_2= z_{(\{1,1\},\{3\})},\quad z_3=z_{(\{1,1\},\{1\})},\quad z_4=z_{(\{2,2\},\{2\})}.\]
We have
\[
 \pi_{\mu}(z_1\cdot z_2\cdot z_3\cdot z_4) = m(z_1,z_2)\cdot m(z_3,z_4) +m(z_1,z_3)\cdot m(z_2,z_4) +m(z_1,z_4)\cdot m(z_2,z_3)=\]
\[z_{(\{1,1,1,2\},\{1,3\})}\cdot z_{(\{1,1,2,2\},\{1,2\})}+z_{(\{1,1,1,2\},\{1,1\})}\cdot z_{(\{1,1,2,2\},\{2,3\})}+z_{(\{1,2,2,2\},\{1,2\})}\cdot z_{(\{1,1,1,1\},\{1,3\})}.
\]
\end{example}

A more ``visual'' way of representing the monomials in $\Sym(\Sym^{d_1}V_1\o\cdots\o\Sym^{d_n}V_n)=K[z_{\a}]$ and the maps $\pi_{\mu}$ is as follows. We identify each $z_{\a}$ with an $1\times n$ block with entries the multisets $\a_i$:
\[
 z_{\a} = \begin{array}{|c|c|c|c|}
 \hline
 \a_1 & \a_2 & \cdots & \a_n \\ \hline
 \end{array}.
\]
We represent a monomial $m = z_{\a^1}\cdots z_{\a^r}$ of degree $r$ as an $r\times n$ block $M$, whose rows correspond to the variables $z_{\a^i}$ in the way described above.
\[
 m \equiv M = \begin{array}{|c|c|c|c|}
 \hline
 \a^1_1 & \a^1_2 & \cdots & \a^1_n \\ \hline
 \a^2_1 & \a^2_2 & \cdots & \a^2_n \\ \hline
 \vdots & \vdots & \ddots & \vdots \\ \hline
 \a^r_1 & \a^r_2 & \cdots & \a^r_n \\ \hline
 \end{array}
\]
Note that the order of the rows is irrelevant, since the $z_{\a^i}$'s commute. The way $\pi_{\mu}$ acts on an $r\times n$ block $M$ is as follows: it partitions in all possible ways the set of rows of $M$ into subsets of sizes equal to the parts of $\mu$, collapses the elements of each subset into a single row, and takes the sum of all blocks obtained in this way. Here by collapsing a set of rows we mean taking the columnwise union of the entries of the rows. More precisely, if $M$ is the $r\times n$ block corresponding to $z_{\a^1}\cdots z_{\a^r}$ and $\mu=(\mu_1,\cdots,\mu_t)$, then
\[\pi_{\mu}M=\sum_{\substack{A\in\mc{P}_{\mu}\\A=\{A_1,\cdots,A_t\}}}
 \begin{array}{|c|c|c|}
  \hline
 \cdots & \bigcup_{i\in A_1}\a_k^i & \cdots\\ \hline
 \cdots & \bigcup_{i\in A_2}\a_k^i & \cdots\\ \hline
 \vdots & \ddots & \vdots \\ \hline
 \cdots & \bigcup_{i\in A_t}\a_k^i & \cdots\\ \hline
 \end{array}.
\]
Note that if two $A_i$'s have the same cardinality, then the variables corresponding to their rows commute, so we can harmlessly interchange them.

\begin{example}\label{ex:pimublock} With these conventions, we can rewrite Example \ref{ex:pimu} as 
\[
\begin{array}{|c|c|}
 \hline
1,2 & 1 \\ \hline
1,1 & 3 \\ \hline
1,1 & 1 \\ \hline
2,2 & 2 \\ \hline
\end{array}
\overset{\pi_{(2,2)}}{\lra}
\begin{array}{|c|c|}
 \hline
1,1,1,2 & 1,3 \\ \hline
1,1,2,2 & 1,2 \\ \hline
\end{array}
+
\begin{array}{|c|c|}
 \hline
1,1,1,2 & 1,1 \\ \hline
1,1,2,2 & 2,3 \\ \hline
\end{array}
+
\begin{array}{|c|c|}
 \hline
1,2,2,2 & 1,2 \\ \hline
1,1,1,1 & 1,3 \\ \hline
\end{array}.
\]
\end{example}

\begin{proposition}[Multi-prolongations, \cite{landsberg}]\label{prop:prolongations} For a positive integer $r$, the polynomials of degree $r$ vanishing on $\s_k(SV_{d_1,\cdots,d_n}(\bb{P}V_1^*\times\cdots\times\bb{P}V_n^*))$ are precisely the elements of $S_{(r)}(S_{(d_1)}V_1\o\cdots\o S_{(d_n)}V_n)$ in the intersection of the kernels of the maps $\pi_{\mu}$, where $\mu$ ranges over all partitions of $r$ with (at most) $k$ parts.
\end{proposition}

\begin{proof}
 Let $X$ denote the Segre-Veronese variety $SV_{d_1,\cdots,d_n}(\bb{P}V_1^*\times\cdots\times\bb{P}V_n^*)$. As in Section \ref{sec:secants}, there exists a map ($s^{\#}$, which we now denote $\pi$)
\[\pi:\Sym(S_{(d_1)}V_1\o\cdots\o S_{(d_n)}V_n)\lra K[X]^{\o k},\]
whose kernel and image coincide with the ideal and homogeneous coordinate ring respectively, of $\s_k(X)$. Using the description of $K[X]$ given in Section \ref{sec:coordringSV}, we obtain that the degree $r$ part of the target of $\pi$ is
\[(K[X]^{\o k})_r = \bigoplus_{\mu_1+\cdots+\mu_k = r}(S_{(\mu_1 d_1)}V_1\o\cdots\o S_{(\mu_1 d_n)}V_n)\o\cdots\o(S_{(\mu_k d_1)}V_1\o\cdots\o S_{(\mu_k d_n)}V_n).\]

The degree $r$ component of $\pi$, which we call $\pi_r$, is then a direct sum of maps $\pi_{\mu}$ as in Definition \ref{def:pimu}, where $\mu$ ranges over partitions of $r$ with at most $k$ parts. The conclusion of the proposition now follows. To see that it's enough to only consider partitions with exactly $k$ parts, note that if $\mu$ has fewer than $k$ parts, and $\widehat{\mu}$ is a partition obtained by subdividing $\mu$ (splitting some of the parts of $\mu$ into smaller pieces), then $\pi_{\mu}$ factors through (up to a multiplicative factor) $\pi_{\widehat{\mu}}$, hence $\rm{ker}(\pi_{\mu})\supset \rm{ker}(\pi_{\widehat{\mu}})$, so the contribution of $\rm{ker}(\pi_{\mu})$ to the intersection of kernels is superfluous.
\end{proof}

\begin{definition}[Multi-prolongations]\label{def:prolongations}
 We write $I_{\mu}(V)=I_{\mu}^{\d}(V)$ for the kernel of the map $\pi_{\mu}(V)$, and $I_{r}(V)=I_r^{\d}(V)$ for the intersection of the kernels of the maps $\pi_{\mu}(V)$ as $\mu$ ranges over partitions of $r$ with $k$ parts. I.e. $I_{r}(V)$ is the degree $r$ part of the ideal of $\s_k(SV_{d_1,\cdots,d_n}(\bb{P}V_1^*\times\cdots\times\bb{P}V_n^*))$.
\end{definition}

Given the description of the ideal of $\s_k(X)$ as the kernel of the $GL(V)$-equivariant map $\pi$, we now proceed to analyze $\pi$ irreducible representation by irreducible representation. That is, we fix a positive integer $r$ and an $n$-partition $\ll=(\ll^1,\cdots,\ll^n)$ of $(rd_1,\cdots,rd_n)$, and we restrict $\pi$ to the $\ll$-parts of its source and target. The map $\pi$ depends functorially on the vector spaces $V_1,\cdots,V_n$, and its kernel and image stabilize from a representation theoretic point of view as the dimensions of the $V_i$'s increase. More precisely, we have the following

\begin{proposition}[Inheritance, \cite{landsberg}]\label{prop:inheritance}
 Fix an $n$-partition $\ll\vdash^n r\cdot(d_1,\cdots,d_n)$. Let $l_j$ denote the number of parts of $\ll^j$, for $j=1,\cdots,n$. Then the multiplicities of $S_{\ll}V$ in the kernel and image respectively of $\pi$ are independent of the dimensions $m_j$ of the $V_j$'s, as long as $m_j\geq l_j$. Moreover, if some $l_j$ is larger than $k$, then $S_{\ll}V$ doesn't occur as a representation in the image of $\pi$.
\end{proposition}

\begin{proof}
 The last statement follows from the representation theoretic description of the coordinate ring of a Segre-Veronese variety, and Pieri's rule: every irreducible representation $S_{\ll}V$ occurring in $K[X]^{\o k}$ must have the property that each $\ll^j$ has at most $k$ parts.

As for the first part, note that $\pi$ is completely determined by what it does on the $\ll$-highest weight vectors, and that the $\ll$-highest weight vector of an irreducible representation $S_{\ll}V$ only depends on the first $l_j$ elements of the basis $\mc{B}_j$, for $j=1,\cdots,n$. 
\end{proof}

We just saw in the previous proposition that (the $\ll$-part of) $\pi$ is essentially insensitive to expanding or shrinking the vector spaces $V_i$, as long as their dimensions remain larger than $l_i$. Also, the last part of the proposition allows us to concentrate on $n$-partitions $\ll$ where each $\ll^i$ has at most $k$ parts. To understand $\pi$, we thus have the freedom to pick the dimensions of the $V_i$'s to be positive integers at least equal to $k$. It might seem natural then to pick these dimensions as small as possible (equal to $k$), and understand the kernel and image of $\pi$ in that situation. However, we choose not to do so, and instead we fix a positive degree $r$ and concentrate our attention on the map $\pi_r$, the degree $r$ part of $\pi$. We assume that
\[\rm{dim}(V_i)=r\cdot d_i,\quad i=1,\cdots,n.\]
The reason for this assumption is that now the $\sl$ zero-weight spaces of the source and target of $\pi_r$ are nonempty and generate the corresponding representations. Therefore $\pi_r$ is determined by its restriction to these zero-weight spaces, which suddenly makes our problem combinatorial: the zero-weight spaces are modules over the Weyl group, which is just the product of symmetric groups $S_{rd_1}\times\cdots\times S_{rd_n}$, allowing us to use the representation theory of the symmetric groups to analyze the map $\pi_r$. We call this reduction the ``generic case'', because the $\sl$ zero-weight subspace of $S_{(r)}(S_{(d_1)}V_1\o\cdots\o S_{(d_n)}V_n)$ is the subspace containing the most generic tensors.

\subsection{The ``generic case''}\label{sec:generic}

\subsubsection{Generic multi-prolongations}\label{subsec:genericmultiprolongations}

We let $\d,\r$ denote the sequences of numbers $(d_1,\cdots,d_n)$ and $r\cdot\d=(rd_1,\cdots,rd_n)$ respectively. We let $S_{\r}$ denote the product of symmetric groups $S_{rd_1}\times\cdots\times S_{rd_n}$, the Weyl group of the Lie algebra $\sl(V)$ of $GL(V)$ (recall that $\rm{dim}(V_j)=m_j=rd_j$ for $j=1,\cdots,n$). 

\begin{definition}\label{def:urd} We denote by $U_{r}^{\d}$ the $\sl(V)$ zero-weight space of the representation $S_{(r)}(S_{(d_1)}V_1\o\cdots\o S_{(d_n)}V_n)$. $U_{r}^{\d}$ has a basis consisting of monomials
\[m = z_{\a^1}\cdots z_{\a^r},\]
where for each $j$, the elements of $\{\a_j^1,\cdots,\a_j^r\}$ form a partition of the set $\{1,\cdots,rd_j\}$, with $|\a_j^i|=d_j$. Alternatively, $U_{r}^{\d}$ has a basis consisting of $r\times n$ blocks $M$, where each column of $M$ yields a partition of the set $\{1,\cdots,rd_j\}$ with $r$ equal parts.
\end{definition}

\begin{example}\label{ex:urd}
 For $n=2$, $d_1=2$, $d_2=1$, $r=4$, a typical element of $U_{r}^{\d}$ is
\[
 M = \begin{array}{|c|c|}
 \hline
1,6 & 1 \\ \hline
2,3 & 4 \\ \hline
4,5 & 2 \\ \hline
7,8 & 3 \\ \hline
\end{array} = 
z_{(\{1,6\},\{1\})}\cdot z_{(\{2,3\},\{4\})}\cdot z_{(\{4,5\},\{2\})}\cdot z_{(\{7,8\},\{3\})} = m.
\]
\end{example}

$S_{\r}$ acts on $U_{r}^{\d}$ by letting its $j$-th factor $S_{rd_j}$ act on the $j$-th columns of the blocks $M$ described above. As an abstract representation, we have
\[U_{r}^{\d}\simeq\rm{Ind}_{(S_{d_1}\times\cdots\times S_{d_n})^{r}\wr S_r}^{S_{\r}}({\bf 1}),\]
where $\wr$ denotes the \defi{wreath product} of $(S_{d_1}\times\cdots\times S_{d_n})^{r}$ with $S_r$, and ${\bf 1}$ denotes the trivial representation (we will say more about this in the following section). The dimension of the space $U_{r}^{\d}$ is
\[
 N = N_{r}^{\d} = \frac{(rd_1)!(rd_2)!\cdots(rd_n)!}{(d_1!d_2!\cdots d_n!)^r\cdot r!}. 
\]

\begin{example}\label{ex:actionSr} Continuing Example \ref{ex:urd}, let $\s=(\s_1,\s_2)\in S_8\times S_4$, with $\s_1=(1,2)(5,3,7)$, $\s_2=(1,4,3)$, in cycle notation. Then
\[\s\cdot M = \begin{array}{|c|c|}
 \hline
2,6 & 4 \\ \hline
1,7 & 3 \\ \hline
4,3 & 2 \\ \hline
5,8 & 1 \\ \hline
\end{array},
\]
or
\[\s\cdot m = z_{(\{2,6\},\{4\})}\cdot z_{(\{1,7\},\{3\})}\cdot z_{(\{4,3\},\{2\})}\cdot z_{(\{5,8\},\{1\})}.\]
\end{example}

\begin{definition}\label{def:upid}
 For a partition $\mu$ written in multiplicative notation $\mu=(\mu_1^{i_1}\cdots\mu_s^{i_s})$ as in Definition \ref{def:pimu}, we define the space $U_{\mu}^{\d}$ to be the $\sl$ zero-weight space of the representation 
\[\bigotimes_{j=1}^s S_{(i_j)}(S_{(\mu_j d_1)}V_1\o\cdots\o S_{(\mu_j d_n)}V_n).\]
Writing $\mu=(\mu_1,\cdots,\mu_t)$ we can realize $U_{\mu}^{\d}$ as the vector space with a basis consisting of $t\times n$ blocks $M$ with the entry in row $i$ and column $j$ consisting of $\mu_i\cdot d_j$ elements from the set $\{1,\cdots,rd_j\}$, in such a way that each column of $M$ represents a partition of $\{1,\cdots,rd_j\}$. As usual, we identify two blocks if they differ by permutations of rows of the same size, i.e. corresponding to equal parts of $\mu$. Note that when $\mu=(1^r)$ we get $U_{\mu}^{\d}=U_{r}^{\d}$, recovering Definition \ref{def:urd}.
\end{definition}

We can now define the generic version of the map $\pi_{\mu}$ from Definition \ref{def:pimu}:

\begin{definition}\label{def:pimugeneric}
 For a partition $\mu\vdash r$ as in Definition \ref{def:pimu}, we define the map
\[\pi_{\mu}:U_{r}^{\d}\lra U_{\mu}^{\d},\]
to be the restriction of the map from Definition \ref{def:pimu} to the $\sl$ zero-weight spaces of the source and target. 
\end{definition}

\begin{example}\label{ex:pimugeneric} The generic analogue of Example \ref{ex:pimublock} is:
\[
\begin{array}{|c|c|}
 \hline
1,6 & 1 \\ \hline
2,3 & 4 \\ \hline
4,5 & 2 \\ \hline
7,8 & 3 \\ \hline
\end{array}
\overset{\pi_{(2,2)}}{\lra}
\begin{array}{|c|c|}
 \hline
1,2,3,6 & 1,4 \\ \hline
4,5,7,8 & 2,3 \\ \hline
\end{array}
+
\begin{array}{|c|c|}
 \hline
1,4,5,6 & 1,2 \\ \hline
2,3,7,8 & 3,4 \\ \hline
\end{array}
+
\begin{array}{|c|c|}
 \hline
1,6,7,8 & 1,3 \\ \hline
2,3,4,5 & 2,4 \\ \hline
\end{array}.
\]
If instead of the partition $(2,2)$ we take $\mu=(2,1,1)=(1^2 2)$, then we have
\[
\begin{array}{|c|c|}
 \hline
1,6 & 1 \\ \hline
2,3 & 4 \\ \hline
4,5 & 2 \\ \hline
7,8 & 3 \\ \hline
\end{array}
\overset{\pi_{(2,1,1)}}{\lra}
\begin{array}{|c|c|}
 \hline
1,2,3,6 & 1,4 \\ \hline
4,5 & 2 \\ \hline
7,8 & 3 \\ \hline
\end{array}
+
\begin{array}{|c|c|}
 \hline
1,4,5,6 & 1,2 \\ \hline
2,3 & 4 \\ \hline
7,8 & 3 \\ \hline
\end{array}
+\]
\[
\begin{array}{|c|c|}
 \hline
1,6,7,8 & 1,3 \\ \hline
2,3 & 4 \\ \hline
4,5 & 2 \\ \hline
\end{array}
+
\begin{array}{|c|c|}
 \hline
2,3,4,5 & 2,4 \\ \hline
1,6 & 1 \\ \hline
7,8 & 3 \\ \hline
\end{array}
+
\begin{array}{|c|c|}
 \hline
2,3,7,8 & 3,4 \\ \hline
1,6 & 1 \\ \hline
4,5 & 2 \\ \hline
\end{array}
+
\begin{array}{|c|c|}
 \hline
4,5,7,8 & 2,3 \\ \hline
1,6 & 1 \\ \hline
2,3 & 4 \\ \hline
\end{array}.
\] 
Note that if we compose $\pi_{(2,1,1)}$ with the multiplication map that collapses together the last two rows of a block in $U_{(2,1,1)}^{(2,1)}$, then we obtain the map $2\cdot\pi_{(2,2)}$.
\end{example}

\begin{definition}[Generic multi-prolongations]\label{def:genericprolongations}
We write $I_{\mu}=I_{\mu}^{\d}$ for the kernel of $\pi_{\mu}$, and $I_r=I_r^{\d}$ for the intersection of the kernels of the maps $\pi_{\mu}$, as $\mu$ ranges over partitions of $r$ with at most (exactly) $k$ parts. We refer to $I_{r}$ as the set of ``generic equations'' for $\s_k(SV_{\d}(\bb{P}V_1^*\o\cdots\o\bb{P}V_n^*))$, or ``generic multi-prolongations'' (see Proposition \ref{prop:prolongations}).
\end{definition}

\subsubsection{Tableaux}\label{subsec:tableaux}

The maps $\pi_{\mu}$, for various partitions $\mu$, are $S_{\r}$-equivariant, so to understand them it suffices to analyze them irreducible representation by irreducible representation. Recall that irreducible $S_{\r}$-representations are classified by $n$-partitions $\ll\vdash^n\r$, so we fix one such. This gives rise to a Young symmetrizer $c_{\ll}$ as explained in Section \ref{sec:reptheory}, and all the data of $\pi_{\mu}$ (concerning the $\ll$-parts of its kernel and image) is contained in its restriction to the $\ll$-highest weight spaces of the source and target, i.e. in the map
\[\pi_{\mu} = \pi_{\mu}(\ll):c_{\ll}\cdot U_{r}^{\d}\lra c_{\ll}\cdot U_{\mu}^{\d}.\]

We now introduce the tableaux formalism that's fundamental for the proof of our main results, giving a combinatorial perspective on the analysis of the kernels and images of the maps $\pi_{\mu}$, which are the main objects we're after.

The representations $U_{\mu}^{\d}$ are spanned by blocks $M$ as in Definition \ref{def:upid}, hence the vector spaces  $c_{\ll}\cdot U_{\mu}^{\d}$ are spanned by elements of the form $c_{\ll}\cdot M$, which we shall represent as $n$-tableaux, according to the following definition.

\begin{definition}\label{def:tableaux}
 Given a partition $\mu=(\mu_1,\cdots,\mu_t)\vdash r$, an $n$-partition $\ll\vdash^n\r$ and a block $M\in U_{\mu}^{\d}$, we associate to the element $c_{\ll}\cdot M\in c_{\ll}\cdot U_{\mu}^{\d}$ the $n$-tableau 
\[T=(T^1,\cdots,T^n)=T^1\o\cdots\o T^n\]
of shape $\ll$, obtained as follows. Suppose that the block $M$ has the set $\a_j^i$ in its $i$-th row and $j$-th column. Then we set equal to $i$ the entries in the boxes of $T^j$ indexed by elements of $\a_j^i$ (recall from Section \ref{sec:reptheory} that the boxes of a tableau are indexed canonically: from left to right and top to bottom). Note that each tableau $T^j$ has entries $1,\cdots,t$, with $i$ appearing exactly $\mu_i\cdot d_j$ times.

Note also that in order to construct the $n$-tableau $T$ we have made a choice of the ordering of the rows of $M$: interchanging rows $i$ and $i'$ when $\mu_i=\mu_{i'}$ should yield the same element $M\in U_{\mu}^{\d}$, therefore we identify the corresponding $n$-tableaux that differ by interchanging the entries equal to $i$ and $i'$. 
\end{definition}

\begin{example}\label{ex:tableaux} We let $n=2,\d=(2,1),r=4,\mu=(2,2)$ as in Example \ref{ex:pimu}, and consider the $2$-partition $\ll=(\ll^1,\ll^2)$, with $\ll^1=(5,3)$, $\ll^2=(2,1,1)$. We have
\[
\xymatrix
{
{c_{\ll}\cdot\begin{array}{|c|c|}
 \hline
1,6 & 1 \\ \hline
2,3 & 4 \\ \hline
4,5 & 2 \\ \hline
7,8 & 3 \\ \hline
\end{array}} \ar@{=}[r]\ar@{=}[d] &  {\Yvcentermath1\young(12233,144)\o\young(13,4,2)} \ar@{=}[d]\\
{c_{\ll}\cdot\begin{array}{|c|c|}
 \hline
2,3 & 4 \\ \hline
7,8 & 3 \\ \hline
1,6 & 1 \\ \hline
4,5 & 2 \\ \hline
\end{array}} \ar@{=}[r] & {\Yvcentermath1\young(31144,322)\o\young(34,2,1)}\\
}
\]
Let's write down the action of the map $\pi_{\mu}$ on the tableaux pictured above
\[\Yvcentermath1\begin{split}
 \pi_{\mu}\left(\young(12233,144)\o\young(13,4,2)\right)=\ &\young(11122,122)\o\young(12,2,1)+\young(12211,122)\o\young(11,2,2)\\
&+\young(12222,111)\o\young(12,1,2).
\end{split}
\]
\end{example}

We collect in the following lemma the basic relations that $n$-tableaux satisfy.

\begin{lemma}\label{lem:relstableaux}
 Fix an $n$ partition $\ll\vdash^n\r$, and let $T$ be an $n$-tableau of shape $\ll$. The following relations hold:
\begin{enumerate}
 \item If $\s$ is a permutation of the entries of $T$ that preserves the set of entries in each column of $T$, then
\[\s(T)=\rm{sgn}(\s)\cdot T.\]
In particular, if $T$ has repeated entries in a column, then $T=0$.

 \item If $\s$ is a permutation of the entries of $T$ that interchanges columns of the same size of some tableau $T^j$, then
\[\s(T)=T.\]

 \item Assume that one of the tableaux of $T$, say $T^j$ has a column $C$ of size $t$ with entries $a_1,a_2,\cdots,a_t$, and that $b$ is an entry of $T^j$ to the right of $C$. Let $\s_i$ denote the transposition that interchanges $a_i$ with $b$. We have
\[T=\sum_{i=1}^t\s_i(T).\]
We write this as
\[\Yvcentermath1\Yboxdim{18pt}\young(\aone b,\vdots,\ai,\vdots,\at)=\sum_{i=1}^t\young(\aone\ai,\vdots,b,\vdots,\at),\]
disregarding the entries of $T$ that don't get perturbed.
\end{enumerate}
\end{lemma}

\begin{proof}
 (1) follows from the fact that if $\s\in C_{\ll}$ is a column permutation, then $b_{\ll}\cdot\s=-b_{\ll}$.

 (2) follows from the fact that if $\s$ permutes columns of the same size, then $\s\in R_{\ll}$ is a permutation that preserves the rows of the canonical $n$-tableau of shape $\ll$ (so in particular $a_{\ll}\cdot\s=a_{\ll}$), and $\s$ commutes with $b_{\ll}$. It follows that
\[c_{\ll}\cdot\s=a_{\ll}\cdot(b_{\ll}\cdot\s)=a_{\ll}\cdot(\s\cdot b_{\ll})=(a_{\ll}\cdot\s)\cdot b_{\ll}=a_{\ll}\cdot b_{\ll}=c_{\ll}.\]

 (3) follows from Corollary \ref{cor:1flattenings} (note the rest of the proof uses the formalism of Section \ref{subsec:genericflattenings} below). Let us assume first that all entries $a_1,\cdots,a_t,b$ are distinct. If $\tilde{T}$ is the $n$-tableau obtained by circling the entries $a_1,\cdots,a_t,b$, then
\[\Yvcentermath1\Yboxdim{18pt}\tilde{T} = \young(\ccaone b,\vdots,\ccai,\vdots,\ccat) - \sum_{i=1}^t\young(\ccaone\ai,\vdots,\ccb,\vdots,\ccat).\]
By skew-symmetry on columns (part (1)), the effect of circling $t$ entries in the same column of a tableau $T$ is precisely multiplying $T$ by $t!$. It follows that we can rewrite the above relation as
\[\Yvcentermath1\Yboxdim{18pt}\tilde{T}=t!\cdot\left(\young(\aone b,\vdots,\ai,\vdots,\at)-\sum_{i=1}^t\young(\aone\ai,\vdots,b,\vdots,\at)\right).\]
By Corollary \ref{cor:1flattenings}, $\tilde{T}=0$, which combined with the above equality yields the desired relation.

Now if $a_1,\cdots,a_t,b$ are not distinct, then either $a_i=a_j$ for some $i\neq j$, or $b=a_i$ for some $i$. If $a_i=a_j$, then $T$ and $\s_k(T)$, $k\neq i,j$, have repeated entries in the column $C$, hence they are zero. Relation (3) becomes then $0=\s_i(T)+\s_j(T)$. But this is true by part (1), because $\s_i(T)$ and $\s_j(T)$ differ by a column transposition.

Assume now that $b=a_i$ for some $i$. Then $\s_j(T)$ has repeated entries in the column $C$ for $j\neq i$, thus relation (3) becomes $T=\s_i(T)$, which is true because $a_i=b$.
\end{proof}

There is one last ingredient that we need to introduce in the generic setting, namely the \defi{generic flattenings}.

\subsubsection{Generic flattenings}\label{subsec:genericflattenings}

\begin{definition}[Generic flattenings]\label{def:genericflattenings} For a decomposition $\d=A+B$, $A=(a_1,\cdots,a_n)$, $B=(b_1,\cdots,b_n)$, (i.e. $d_i=a_i+b_i$ for $i=1,\cdots,n$), we write $F_{A,B}^{k,r}$ for the span of \defi{$k\times k$ minors of generic $(A,B)$-flattenings}. This is the subspace of $U_{r}^{\d}$ spanned by expressions of the form
\[[\a^1,\cdots,\a^k|\b^1,\cdots,\b^k]\cdot z_{\c^{k+1}}\cdots z_{\c^r},\]
where $[\a^1,\cdots,\a^k|\b^1,\cdots,\b^k]=\rm{det}(z_{\a^i\cup\b^j})$, $\a^i=(\a^i_1,\cdots,\a^i_n)$, $\b^i=(\b^i_1,\cdots,\b^i_n)$, $\c^i=(\c^i_1,\cdots,\c^i_n)$, with $|\a^i_j|=a_j$, $|\b^i_j|=b_j$ and $|\c^i_j|=d_j$, and such that for fixed $j$, the sets  $\a^i_j,\b^i_j,\c^i_j$ form a partition of the set $\{1,\cdots,rd_j\}$.
\end{definition}

\begin{example}\label{ex:genericflattenings} Take $n=2$, $\d=(2,1)$ and $r=4$, as usual. Take $A=(1,1)$, $B=(1,0)$ and $k=3$. A typical element of $F_{A,B}^{3,4}$ looks like
\[D=[(\{1\},\{1\}),(\{3\},\{4\}),(\{7\},\{3\})|(\{6\},\{\}),(\{2\},\{\}),(\{8\},\{\})]\cdot z_{(\{4,5\},\{2\})}=\]
\[\rm{det}\left[ \begin{array}{ccc}
z_{(\{1,6\},\{1\})} & z_{(\{1,2\},\{1\})} & z_{(\{1,8\},\{1\})} \\
z_{(\{3,6\},\{4\})} & z_{(\{3,2\},\{4\})} & z_{(\{3,8\},\{4\})} \\
z_{(\{7,6\},\{3\})} & z_{(\{7,2\},\{3\})} & z_{(\{7,8\},\{3\})} 
\end{array} \right]\cdot z_{(\{4,5\},\{2\})}.\]
Expanding the determinant, we obtain
\[D=\begin{array}{|c|c|}
 \hline
1,6 & 1 \\ \hline
3,2 & 4 \\ \hline
7,8 & 3 \\ \hline
4,5 & 2 \\ \hline
\end{array}-
\begin{array}{|c|c|}
 \hline
1,2 & 1 \\ \hline
3,6 & 4 \\ \hline
7,8 & 3 \\ \hline
4,5 & 2 \\ \hline
\end{array}-
\begin{array}{|c|c|}
 \hline
1,8 & 1 \\ \hline
3,2 & 4 \\ \hline
7,6 & 3 \\ \hline
4,5 & 2 \\ \hline
\end{array}-
\begin{array}{|c|c|}
 \hline
1,6 & 1 \\ \hline
3,8 & 4 \\ \hline
7,2 & 3 \\ \hline
4,5 & 2 \\ \hline
\end{array}+
\begin{array}{|c|c|}
 \hline
1,8 & 1 \\ \hline
3,6 & 4 \\ \hline
7,2 & 3 \\ \hline
4,5 & 2 \\ \hline
\end{array}+
\begin{array}{|c|c|}
 \hline
1,2 & 1 \\ \hline
3,8 & 4 \\ \hline
7,6 & 3 \\ \hline
4,5 & 2 \\ \hline
\end{array}.
\]
Notice that all the blocks in the above expansion coincide, except in the entries $2,6,8$ that get permuted in all possible ways. Let's multiply now $D$ with the Young symmetrizer $c_{\ll}$ for $\ll=(\ll^1,\ll^2)$, $\ll^1=(5,3)$ and $\ll^2=(2,1,1)$. We get
\[\Yvcentermath1 c_{\ll}\cdot D = \young(12244,133)\o\young(14,3,2)-\young(11244,233)\o\young(14,3,2)-\young(12244,331)\o\young(14,3,2)\]
\[\Yvcentermath1 -\young(13244,132)\o\young(14,3,2)+\young(13244,231)\o\young(14,3,2)+\young(11244,332)\o\young(14,3,2).\]
Note that all the $2$-tableaux in the previous expression coincide, except in the $2$-nd, $6$-th and $8$-th box of their first tableau, which get permuted in all possible ways. We represent $c_{\ll}\cdot D$ by a $2$-tableau with the entries in boxes $2$, $6$ and $8$ of its first tableau circled (see also Definition \ref{def:circledtableaux} below):
\[\Yvcentermath1 c_{\ll}\cdot D= \young(1\cctwo244,\ccone3\ccthree)\o\young(14,3,2).\]
To reformulate this one last time, we write
\[\Yvcentermath1\young(1\cctwo244,\ccone3\ccthree)\o\young(14,3,2)=\sum_{\s\in S_3}\rm{sgn}(\s)\cdot \s\left(\young(12244,133)\o\young(14,3,2)\right),\]
where $S_{3}=S_{\{1,2,3\}}$ is the symmetric group on the circled entries.
\end{example}

\begin{definition}\label{def:circledtableaux}
 Let $A,B$ and $F_{A,B}^{k,r}$ as in Definition \ref{def:genericflattenings}, let
\[D = [\a^1,\cdots,\a^k|\b^1,\cdots,\b^k]\cdot z_{\c^{k+1}}\cdots z_{\c^r}\in F_{A,B}^{k,r},\]
and let $\ll\vdash^n\r=(rd_1,\cdots,rd_n)$. We let $\c^i=\a^i\cup\b^i$ for $i=1,\cdots,k$, and consider $T=c_{\ll}\cdot m$ the $n$-tableau corresponding to the monomial
\[m=z_{\c^1}\cdots z_{\c^r}.\]
We represent $c_{\ll}\cdot D\in\rm{hwt}_{\ll}(F_{A,B}^{k,r})$ as the $n$-tableau $T$ with the entries in the boxes corresponding to the elements of $\a^1,\cdots,\a^k$ circled. Alternatively, we can circle the entries in the boxes corresponding to the elements of $\b^1,\cdots,\b^k$.
\end{definition}

It follows that a spanning set for $\rm{hwt}_{\ll}(F_{A,B}^{k,r})$ can be obtained as follows: take all the subsets $\mc{C}\subset\{1,\cdots,r\}$ of size $k$, and consider all the $n$-tableaux $T$ with $a_j$ (alternatively $b_j$) of each of the elements of $\mc{C}$ circled in $T^j$. Of course, because of the symmetry of the alphabet $\{1,\cdots,r\}$, it's enough to only consider $\mc{C}=\{1,\cdots,k\}$, so that the only entries we ever circle are $1,2,\cdots,k$.

Continuing with Example \ref{ex:genericflattenings}, we have
\[\Yvcentermath1 c_{\ll}\cdot D = \young(1\cctwo244,\ccone3\ccthree)\o\young(14,3,2) = \young(\ccone2\cctwo44,1\ccthree3)\o\young(\ccone4,\ccthree,\cctwo).\]

Our goal is to reduce the statement of Theorem \ref{thm:main} to an equivalent statement that holds in the generic setting, and thus transform our problem into a combinatorial one. More precisely, we would like to say that the space of generic flattenings coincides with the intersection of the kernels of the (generic) maps $\pi_{\mu}$, and that this is enough to conclude the same about the nongeneric case. One issue that arises is that we don't know at this point (although it seems very tempting to assert) that the zero-weight space of the space of flattenings coincides with the space of generic flattenings. Section \ref{sec:polarization} will show how to take care of this issue, and how to reduce all our questions to the generic setting. 

\subsubsection{$1$-flattenings}\label{subsec:1flattenings}

In this section we focus on the space of generic $1$-flattenings, $F_1=F_1^{k,r}$, defined as the subspace of $U_r^{\d}$ given by
\[F_1^{k,r}=\sum_{\substack{A+B=\d\\|A|=1}}F_{A,B}^{k,r}.\]
We shall see that $F_1$ has a very simple representation theoretic description, which by the results of the next section will carry over to the nongeneric case.

\begin{proposition}\label{prop:1flattenings}
 With the above notations, we have
\[F_1=\bigoplus_{\substack{\ll\vdash^n\r\\ \ll_k\neq 0}}(U_r^{\d})_{\ll},\]
where $(U_r^{\d})_{\ll}$ denotes the $\ll$-part of the representation $U_r^{\d}$, and $\ll_k\neq 0$ means $\ll^j_k\neq 0$ for some $j=1,\cdots,n$, i.e. some partition $\ll^j$ has at least $k$ parts.
\end{proposition}

\begin{proof}
 We divide the proof into two parts:

a) If $\ll\vdash^n\r$ is an $n$-partition with some $\ll^j$ having at least $k$ parts, and $T$ is an $n$-tableau of shape $\ll$, then $T\in F_1$.

b) If $\ll\vdash^n\r$ is an $n$-partition with all $\ll^j$ having less than $k$ parts, then $c_{\ll}\cdot F_1=0$.

Let us start by proving part a). We assume that $\ll^j$ has at least $k$ parts and consider $T$ an $n$-tableau of shape $\ll$. If $T^j$ has repeated entries in its first column, then $T=0$. Otherwise, we may assume that the first column of $T^j$ has entries $1,2,\cdots,t$ in this order, where $t$ is the number of parts of $\ll^j$, $t\geq k$. We consider the $n$-tableau $\tilde{T}$ obtained from $T$ by circling the entries $1,2,\cdots,k$ in the first column of $T^j$. We have
\[\Yvcentermath1\Yboxdim{22pt} \tilde{T}=T^1\o\cdots\o\young(\ccone\cdots,\cctwo\cdots,\vdots\vdots,\cck\cdots,\kplone\cdots,\vdots)\o\cdots\o T^n,\]
i.e.
\[\tilde{T}=\sum_{\s\in S_k}\rm{sgn}(\s)\cdot\s(T),\]
where $S_k$ denotes the symmetric group on the circled entries. Since $\s(T)$ differs from $T$ by the column permutation $\s$, it follows by the skew-symmetry of tableaux that
\[\s(T)=\rm{sgn}(\s)\cdot T.\]
This shows that
\[\tilde{T}=k!\cdot T\Longleftrightarrow T=\frac{1}{k!}\cdot\tilde{T}\in F_1,\]
proving a).

To prove b), let
\[D = [\a^1,\cdots,\a^k|\b^1,\cdots,\b^k]\cdot z_{\c^{k+1}}\cdots z_{\c^r}\in F_{A,B}^{k,r},\]
for some $A+B=\d$ with $|A|=1$. We have that $b_{\ll}\cdot D$ is a linear combination of terms that look like $D$, so in order to prove that $c_{\ll}=a_{\ll}\cdot b_{\ll}$ annihilates $D$, it suffices to show that $a_{\ll}\cdot D=0$.

We have $A=(a_1,\cdots,a_n)$ with $a_j=1$ for some $j$ and $a_i=0$ for $i\neq j$. We can thus think of each of $\a^1,\cdots,\a^k$ as specifying a box in the partition $\ll^j$. Since $\ll^j$ has less than $k$ parts, it means that two of these boxes, say $p$ and $q$, lie in the same same row of $\ll^j$. Let $\s=(p,q)$ be the transposition of the two boxes. $\s$ is an element in the group $R_{\ll}$ of permutations that preserve the rows of the canonical $n$-tableau of shape $\ll$ (Section \ref{sec:reptheory}), which means that $a_{\ll}\cdot\s = a_{\ll}$. However,
\[\s\cdot [\a^1,\cdots,\a^p,\cdots,\a^q,\cdots,\a^k|\b^1,\cdots,\b^k] = [\a^1,\cdots,\a^q,\cdots,\a^p,\cdots,\a^k|\b^1,\cdots,\b^k] \]
\[= - [\a^1,\cdots,\a^p,\cdots,\a^q,\cdots,\a^k|\b^1,\cdots,\b^k],\]
since interchanging two rows/columns of a matrix changes the sign of its determinant. We get
\[a_{\ll}\cdot D = (a_{\ll}\cdot\s)\cdot D = a_{\ll}\cdot(\s\cdot D)=a_{\ll}\cdot (-D) = - a_{\ll}\cdot D,\]
hence $a_{\ll}\cdot D=0$, as desired.
\end{proof}

\begin{remark}
 The nongeneric $1$-flattenings give the equations of the so-called \defi{subspace varieties} (see \cite{landsberg} or \cite[Prop.~7.1.2]{weyman}), and in fact this statement is essentially equivalent to our Proposition \ref{prop:1flattenings} via the results of the next section, namely Proposition \ref{prop:genericspecial}.
\end{remark}

\begin{corollary}\label{cor:1flattenings}
 Let $\mc{C}\subset\{1,\cdots,r\}$ be a subset of size $k$. If $\ll$ is an $n$-partition with each $\ll^j$ having less than $k$ parts, and $\tilde{T}$ is an $n$-tableau of shape $\ll$, with one of each entries of $\mc{C}$ in $\tilde{T}^j$ circled, then $\tilde{T}=0$. More generally, with no assumptions on $\ll$, if the circled entries in $\tilde{T}^j$ all lie in columns of size less than $k$, then $\tilde{T}=0$.
\end{corollary}

\begin{proof}
 The first part follows directly from Proposition \ref{prop:1flattenings}, since $\tilde{T}$ is a $1$-flattening, and the space of $1$-flattenings doesn't have nonzero $\ll$-parts when $\ll$ is such that each of its partitions have less than $k$ parts.

For the more general statement, we can apply the argument for part b) of the proof of the previous proposition. If 
\[D = [\a^1,\cdots,\a^k|\b^1,\cdots,\b^k]\cdot z_{\c^{k+1}}\cdots z_{\c^r}\in F_1\]
is such that each $\a^i$ corresponds to a box of $\tilde{T}^j$ situated in a column of size less than $k$, then since column permutations don't change the columns of the boxes corresponding to the $\a^i$'s, it follows that $b_{\ll}\cdot D$ is a combination of expressions $D'$ with the same properties as $D$. To show that $c_{\ll}\cdot D=0$ is thus enough to prove that $a_{\ll}\cdot D=0$. The proof of this statement is identical to the one in the preceding proposition.
\end{proof}

Many of the classical results on the representation theory connected to secant varieties of Segre--Veronese varieties can be recovered from the generic perspective. For some of them, including the Cauchy formula or Strassen's equations, and their generalization by Landsberg and Manivel \cite{lan-man-strassen}, the reader may consult \cite[Chapter 5]{raicu}.

\subsection{Polarization and specialization}\label{sec:polarization}

In this section $V_1,\cdots,V_n$ are again vector spaces of arbitrary dimensions, $\rm{dim}(V_j)=m_j$, $j=1,\cdots,n$. Let $\r=(r_1,\cdots,r_n)$ be a sequence of positive integers, and let 
\[W=V_1^{\o r_1}\o\cdots\o V_n^{\o r_n}.\]
Let $S_{\r}$ denote the product of symmetric groups $S_{r_1}\times\cdots\times S_{r_n}$, and let $G\subset S_{\r}$ be a subgroup. Consider the natural (right) action of $S_{\r}$ on $W$ obtained by letting $S_{r_i}$ act by permuting the factors of $V_i^{\o r_i}$. More precisely, we write the pure tensors in $W$ as
\[v=\bigotimes_{i,j}v_{ij},\quad\rm{with}\ v_{ij}\in V_j,\ j=1,\cdots,n,\ i=1,\cdots,r_j,\]
and for an element $\s=(\s^1,\cdots,\s^n)\in S_{\r}$, we let
\[v*\s = \bigotimes_{i,j}v_{\s^j(i)j}.\]
This action commutes with the (left) action of $GL(V)$ on $W$, and restricts to an action of $G$ on $W$. It follows that $W^G$ is a $GL(V)$-subrepresentation of $W$.

\begin{proposition}\label{prop:polarization} Continuing with the above notation, we let $U=W^G$ and $U'=\rm{Ind}_G^{S_{\r}}({\bf 1})$. Let $\ll\vdash^n\r$ be an $n$-partition with $\ll^j$ having at most $m_j$ parts. The multiplicity of $S_{\ll}V$ in $U$ is the same with that of $[\ll]$ in $U'$.

Moreover, there exist \defi{polarization} and \defi{specialization} maps
\[P_{\ll}:\rm{wt}_{\ll}(U)\lra U',\quad Q_{\ll}:U'\lra\rm{wt}_{\ll}(U),\]
with the following properties:
\begin{enumerate}

\item $Q_{\ll}$ is surjective.

\item $P_{\ll}$ is a section of $Q_{\ll}$.

\item $P_{\ll}$ and $Q_{\ll}$ restrict to maps between $\rm{hwt}_{\ll}(U)$ and $\rm{hwt}_{\ll}(U')$ which are inverse to each other.

\end{enumerate}
\end{proposition}

\begin{proof}
 The first part is a consequence of Schur-Weyl duality (Lemma \ref{lem:SchurWeyl}) and Frobenius reciprocity (Lemma \ref{lem:Frobenius}). We start with the identification
\[U=W^G=\Hom_G({\bf 1},\rm{Res}_G^{S_{\r}}(W)).\]
Using Schur-Weyl duality we get
\[W=V_1^{\o r_1}\o\cdots\o V_n^{\o r_n}=\bigoplus_{\ll\vdash^n\r}[\ll]\o S_{\ll}V,\]
therefore the previous equality becomes
\[U = \bigoplus_{\ll\vdash^n\r}\Hom_G({\bf 1},\rm{Res}_G^{S_{\r}}([\ll]))\o S_{\ll}V.\]
Frobenius reciprocity now yields
\[\Hom_G({\bf 1},\rm{Res}_G^{S_{\r}}([\ll]))=\Hom_{S_{\r}}(\rm{Ind}_G^{S_{\r}}({\bf 1}),[\ll])=\Hom_{S_{\r}}(U',[\ll]).\]
We get
\[U = \bigoplus_{\ll\vdash^n\r}\Hom_{S_{\r}}(U',[\ll])\o S_{\ll}V,\]
hence the multiplicity of $S_{\ll}V$ in $U$ coincides with that of $[\ll]$ in $U'$, as long as $S_{\ll}V\neq 0$, i.e. as long as $m_j$ is at least as large as the number of parts of the partition $\ll^j$.

It follows that the vector spaces $\rm{hwt}_{\ll}(U)$ and $\rm{hwt}_{\ll}(U')$ have the same dimension, equal to the multiplicity of $S_{\ll}V$ and $[\ll]$ in $U$ and $U'$ respectively. We next construct explicit maps $P_{\ll}$, $Q_{\ll}$ inducing isomorphisms of vector spaces between the two spaces.

We identify an element $\s=(\s^1,\cdots,\s^n)\in S_{\r}$ with the ``tensor''
\[\bigotimes_{i,j}\s^j(i),\]
and consider the (regular) representation of $S_{\r}$ on the vector space $R$ with basis consisting of the tensors $\s$ for $\s\in S_{\r}$. The left action of $S_{\r}$ on $R$ is given by
\[\s\cdot\bigotimes_{i,j} a_{ij}=\bigotimes_{i,j}\s^j(a_{ij}),\]
while the right action is given by
\[\bigotimes_{i,j} a_{ij} * \s = \bigotimes_{i,j}a_{\s^j(i)j}.\]
We consider the vector space map $Q_{\ll}:R\to W$ given by
\[\bigotimes_{i,j}a_{ij}\lra\bigotimes_{i,j} g_j(a_{ij}),\]
where $g_j:\{1,\cdots,r_j\}\to\mc{B}_j$ is the map sending $a$ to $x_{ij}$ if the $a$-th box of $\ll^j$ is contained in the $i$-th row of $\ll^j$ (or equivalently if $\ll^j_1+\cdots+\ll^j_{i-1}<a\leq\ll^j_1+\cdots+\ll^j_i$). The image of $Q_{\ll}$ is $\rm{wt}_{\ll}(W)$. It is clear that if $a=\bigotimes_{i,j}a_{ij}$ and $b=\bigotimes_{i,j}b_{ij}$, then $Q_{\ll}(a)=Q_{\ll}(b)$ if and only if $a=\s\cdot b$ for $\s\in S_{\r}$ a permutation that preserves the rows of the canonical $n$-tableau of shape $\ll$. It follows that we can define $P_{\ll}:\rm{wt}_{\ll}(W)\to R$ by
\[P_{\ll}(Q_{\ll}(a))=\frac{1}{\ll!}a_{\ll}\cdot a,\]
where $a_{\ll}$ is the row symmetrizer defined in Section \ref{sec:reptheory}, hence $P_{\ll}$ is a section of $Q_{\ll}$.

Notice that $P_{\ll}$ and $Q_{\ll}$ are maps of right $S_{\r}$-modules, i.e. they respect the $*$-action of $S_{\r}$ on $R$ and $\rm{wt}_{\ll}(W)$ respectively.

Let us prove now that $P_{\ll}$ and $Q_{\ll}$ restrict to inverse isomorphisms between $\rm{hwt}_{\ll}(R)=c_{\ll}\cdot R$ (recall from Section \ref{sec:reptheory} that $c_{\ll}$ denotes the Young symmetrizer corresponding to $\ll$) and $\rm{hwt}_{\ll}(W)$. The two spaces certainly have the same dimension (take $G=\{e\}$ to be the trivial subgroup of $S_{\r}$ and apply the first part of the proposition), so it's enough to prove that for $a'\in\rm{hwt}_{\ll}(R)$

a) $Q_{\ll}(a')\in\rm{hwt}_{\ll}(W)$, and

b) $P_{\ll}(Q_{\ll}(a'))=a'$.

To see why part b) is true, note that
\[P_{\ll}(Q_{\ll}(a_{\ll}\cdot a))=\frac{1}{\ll!}\cdot a_{\ll}^2\cdot a=a_{\ll}\cdot a,\]
i.e. $P_{\ll}\circ Q_{\ll}$ fixes $a_{\ll}\cdot R$. Since $\rm{hwt}_{\ll}(R)=c_{\ll}\cdot R\subset a_{\ll}\cdot R$, it follows that $P_{\ll}(Q_{\ll}(a'))=a'$.
To prove a) we need to show that $Q_{\ll}(a')$ is fixed by the Borel (recall the definition of the Borel subgroup from \ref{sec:reptheory}). It's enough to do this when
\[a'=c_{\ll}\cdot a,\quad a=\bigotimes_{i,j}a_{ij}.\]
The pure tensor $a$ corresponds to an element $\s\in S_{\r}$, so we can write $a=e*\s$, where
\[e=\bigotimes_{i,j}e_{ij}\]
is the ``identity'' tensor, $e_{ij}=i$ for all $i,j$. It follows that
\[Q_{\ll}(a')=Q_{\ll}(c_{\ll}\cdot a)=Q_{\ll}(a_{\ll}\cdot b_{\ll}\cdot e * \s)=\ll!\cdot Q_{\ll}(b_{\ll}\cdot e)*\s.\]
Since the $*$ action commutes with the action of the Borel, it is then enough to prove that $Q_{\ll}(b_{\ll}\cdot e)$ is fixed by the Borel. But this is a direct computation:
\[Q_{\ll}(b_{\ll}\cdot e)=\bigotimes_{i,j} x_{1j}\wedge\cdots\wedge x_{(\ll^j)'_i j},\]
where $(\ll^j)'$ denotes the conjugate partition of $\ll^j$, so that in fact $(\ll^j)'_i$ denotes the number of entries in the $i$-th column of $\ll^j$. In any case, it is clear from the formula of $Q_{\ll}(b_{\ll}\cdot e)$ that it is invariant under the Borel, proving the claim that $P_{\ll}$ and $Q_{\ll}$ restrict to inverse isomorphisms between $\rm{hwt}_{\ll}(W)$ and $\rm{hwt}_{\ll}(R)$.

To finish the proof of the proposition, it suffices to notice that, by Remark \ref{rem:invariants}, we have the identities
\[U=W^G=W*s\quad\rm{and}\quad U'=\rm{Ind}_G^{S_{\r}}({\bf 1})=R*s,\]
where
\[s=\sum_{g\in G}g.\]
Now since $P_{\ll}$, $Q_{\ll}$ respect the $*$ action, it follows that they restrict to inverse isomorphisms between
\[\rm{hwt}_{\ll}(W)*s=\rm{hwt}_{\ll}(W*s)=\rm{hwt}_{\ll}(U)\] 
and 
\[\rm{hwt}_{\ll}(R)*s=\rm{hwt}_{\ll}(R*s)=\rm{hwt}_{\ll}(U'),\]
proving the last part of the proposition.
\end{proof}

We shall apply Proposition \ref{prop:polarization} with $\r=(rd_1,\cdots,rd_n)$ and 
\[U=U_r^{\d}(V)=S_{(r)}(S_{(d_1)}V_1\o\cdots\o S_{(d_n)}V_n),\]
or more generally
\[U=U_{\mu}^{\d}(V)=\bigotimes_{j=1}^s S_{(i_j)}(S_{(\mu_j d_1)}V_1\o\cdots\o S_{(\mu_j d_n)}V_n),\]
the source and target respectively of the map $\pi_{\mu}$ in Definition \ref{def:pimu}. $W$ is now the representation
\[W=V_1^{\o rd_1}\o\cdots\o V_n^{\o rd_n}.\]

We start with $U=U_r^{\d}(V)$. We have $U=W^G$ where
\[G=(S_{d_1}\times\cdots\times S_{d_n})^r\wr S_r\]
is the wreath product between $(S_{d_1}\times\cdots\times S_{d_n})^r$ and $S_r$. Recall that for a group $H$ and positive integer $r$, the \defi{wreath product} $H^r\wr S_r$ of $H^r$ with the symmetric group $S_r$ is just the semidirect product
\[H^r\rtimes S_r,\]
where $S_r$ acts on $H^r$ by permuting the $r$ copies of $H$. We can thus identify an element $\s\in G$ with a collection
\[\s=((\s_j^k)_{\substack{j=1,\cdots,n\\k=1,\cdots,r}},\t),\]
where
\[\s_j^k\in S_{d_j},\quad \t\in S_r.\]
We need to say how we regard $G$ as a subgroup of $S_{\r}$. First of all, we think of $S_{\r}=S_{rd_1}\times\cdots\times S_{rd_n}$ as a product of symmetric groups, where $S_{rd_j}$ acts on the set $\mc{D}_j=\{1,\cdots,rd_j\}$. Then we think of an element $\s\in G$ as an element of $S_{\r}$ by letting $\s_j^k$ act as a permutation of
\[\{(\t(k)-1)\cdot d_j+1,\cdots,\t(k)\cdot d_j\}\subset\mc{D}_j.\]
For example, when $d_1=\cdots=d_n=1$, $G$ is just the group $S_r$, diagonally embedded in $S_r^n$. With this $G$, we let $U'=\rm{Ind}_G^{S_{\r}}({\bf 1})$.

One can now see why the representation $U_{r}^{\d}$, as defined in the previous section, can be identified with $U'$. Recall that $U_r^{\d}$ was defined as a space of $r\times n$ blocks with certain identifications. Consider the block
\[M=\begin{array}{|c|c|c|c|}
 \hline
 \{1,\cdots,d_1\} & \{1,\cdots,d_2\} & \cdots & \{1,\cdots,d_n\} \\ \hline
 \{d_1+1,\cdots,2d_1\} & \{d_2+1,\cdots,2d_2\} & \cdots & \{d_n+1,\cdots,2d_n\} \\ \hline
 \vdots & \vdots & \ddots & \vdots \\ \hline
 \{(r-1)d_1+1,\cdots,rd_1\} & \{(r-1)d_2+1,\cdots,rd_2\} & \cdots & \{(r-1)d_n+1,\cdots,rd_n\} \\ \hline
 \end{array}.
\]
$G$ acts trivially on $M$ (because each $\s_j^k$ does, and because the effect of $\t$ is just permuting the rows of $M$), and all the other blocks are obtained from $M$ by the action of some element of $S_{\r}$. One should think of the span of $M$ thus as the trivial representation $\bf 1$ of $G$ that's induced to $S_{\r}$. 

It is probably best to forget at this point that $U'$ was the zero-weight space of a certain representation, and just think of it abstractly as the induced representation
\[\rm{Ind}_G^{S_{\r}}({\bf 1}),\]
with its realization as a space of blocks. An important point to notice now is that for any decomposition $\d=A+B$, and any $k,r$ we have
\[P_{\ll}(\rm{wt}_{\ll}(F_{A,B}^{k,r}(V)))\subset F_{A,B}^{k,r},\]
and
\[Q_{\ll}(F_{A,B}^{k,r})\subset \rm{wt}_{\ll}(F_{A,B}^{k,r}(V)),\]
where $F_{A,B}^{k,r}$ (Definition \ref{def:genericflattenings}) is the generic version of $F_{A,B}^{k,r}(V)$ (Definition \ref{def:flattenings}). This means that on the corresponding $\ll$-highest weight spaces, $P_{\ll}$ and $Q_{\ll}$ restrict to isomorphisms
\[\rm{hwt}_{\ll}(F_{A,B}^{k,r}(V))\simeq\rm{hwt}_{\ll}(F_{A,B}^{k,r}).\]

\begin{example}
 Here's an example of \defi{specialization}, that involves blocks we're already familiar with. Let $n=2$, $d_1=2$, $d_2=1$, $r=4$, $\ll^1=(5,3)$, $\ll^2=(2,1,1)$. The specialization map $Q_{\ll}$ sends
\[
M = \begin{array}{|c|c|}
 \hline
1,6 & 1 \\ \hline
2,3 & 4 \\ \hline
4,5 & 2 \\ \hline
7,8 & 3 \\ \hline
\end{array}\ \overset{Q_{\ll}}{\lra}\ 
\begin{array}{|c|c|}
 \hline
1,2 & 1 \\ \hline
1,1 & 3 \\ \hline
1,1 & 1 \\ \hline
2,2 & 2 \\ \hline
\end{array} = M'.
\]
$Q_{\ll}$ sends $1,2,3,4,5$ from the first column of $M$ to $1$, because boxes $1,2,3,4,5$ of $\ll^1$ lie in the first row of $\ll^1$, and it sends $6,7,8$ to $2$ because boxes $6,7,8$ of $\ll^1$ lie in its second row. A similar description holds for the second column of $M$ and $\ll^2$.

Although we won't write down explicitly $P_{\ll}(M')$ in this example (see the example below for a concrete illustration of the action of $P_{\ll}$), we will just mention that $P_{\ll}(M')$ is the average of the blocks that specialize to $M'$ via the specialization map $Q_{\ll}$. Of course, $M$ is one such block, but there are many more others.
\end{example}

\begin{example} Let $n=3$, $d_1=d_2=d_3=1$ and $\ll^1=\ll^2=\ll^3=(2,1)$. If \[m=z_{(\{1\},\{1\},\{2\})}z_{(\{2\},\{3\},\{1\})}z_{(\{3\},\{2\},\{3\})}\in U',\ \rm{then}\]
\[Q_{\ll}(m)=z_{(\{1\},\{1\},\{1\})}z_{(\{1\},\{2\},\{1\})}z_{(\{2\},\{1\},\{2\})}\in U\ \rm{and}\]
\[
\begin{split}
P_{\ll}(Q_{\ll}(m))=\frac{1}{8}(&z_{(\{1\},\{1\},\{2\})}z_{(\{2\},\{3\},\{1\})}z_{(\{3\},\{2\},\{3\})}+z_{(\{2\},\{1\},\{2\})}z_{(\{1\},\{3\},\{1\})}z_{(\{3\},\{2\},\{3\})}\\
+&z_{(\{1\},\{1\},\{1\})}z_{(\{2\},\{3\},\{2\})}z_{(\{3\},\{2\},\{3\})}+z_{(\{2\},\{1\},\{1\})}z_{(\{1\},\{3\},\{2\})}z_{(\{3\},\{2\},\{3\})}\\
+&z_{(\{1\},\{2\},\{2\})}z_{(\{2\},\{3\},\{1\})}z_{(\{3\},\{1\},\{3\})}+z_{(\{2\},\{2\},\{2\})}z_{(\{1\},\{3\},\{1\})}z_{(\{3\},\{1\},\{3\})}\\
+&z_{(\{1\},\{2\},\{1\})}z_{(\{2\},\{3\},\{2\})}z_{(\{3\},\{1\},\{3\})}+z_{(\{2\},\{2\},\{1\})}z_{(\{1\},\{3\},\{2\})}z_{(\{3\},\{1\},\{3\})}).
\end{split}
\]
\end{example}

When $U=U_{\mu}^{\d}(V)$, with $\mu=(\mu_1^{i_1}\cdots\mu_s^{i_s})$, we get $U=W^G$, where
\[G=\bigtimes_{j=1}^s \left((S_{\mu_j d_1}\times\cdots\times S_{\mu_j d_n})^{i_j}\wr S_{i_j}\right).\]
It follows that $U'=\rm{Ind}_G^{S_{\r}}({\bf 1})=U_{\mu}^{\d}$ with the realization as a space of blocks explained in the preceding section.

We note that the maps $\pi_{\mu}$ and $\pi_{\mu}(V)$ commute with the polarization and specialization maps $P_{\ll}$, $Q_{\ll}$, i.e. we have a commutative diagram
\begin{equation}\label{eq:pimuspecialization}
 \xymatrix{
 U_r^{\d} \ar@<1ex>[rrr]^{Q_{\ll}} \ar[d]_{\pi_{\mu}} & & & {\rm{wt}_{\ll}(U_r^{\d}(V)) \ar@<1ex>[lll]^{P_{\ll}} \ar[d]^{\pi_{\mu}(V)}}\\
 U_{\mu}^{\d} \ar@<1ex>[rrr]^{Q_{\ll}} & & & {\rm{wt}_{\ll}(U_{\mu}^{\d}(V)) \ar@<1ex>[lll]^{P_{\ll}}}\\
}
\end{equation}

\begin{example}
 Let $\d=(2,1)$, $r=4$, $\mu=(2,2)$, $\ll^1=(5,3)$, $\ll^2=(2,1,1)$. We only illustrate the specialization map $Q_{\ll}$, with the above diagram transposed:
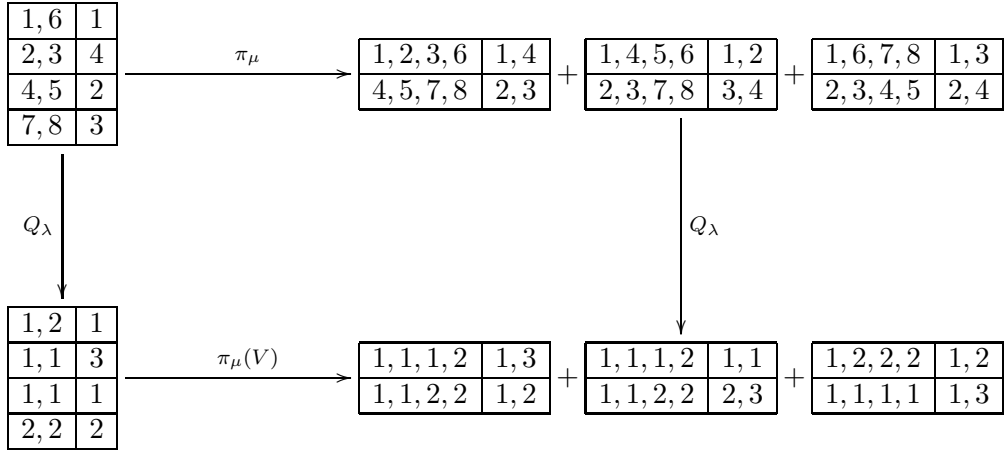
\[
 \xymatrix{
 {\begin{array}{|c|c|}
 \hline
1,6 & 1 \\ \hline
2,3 & 4 \\ \hline
4,5 & 2 \\ \hline
7,8 & 3 \\ \hline
\end{array}} \ar[rrr]^(.3){\pi_{\mu}} \ar[dd]_{Q_{\ll}} & & & {\begin{array}{|c|c|}
 \hline
1,2,3,6 & 1,4 \\ \hline
4,5,7,8 & 2,3 \\ \hline
\end{array}
+
\begin{array}{|c|c|}
 \hline
1,4,5,6 & 1,2 \\ \hline
2,3,7,8 & 3,4 \\ \hline
\end{array}
+
\begin{array}{|c|c|}
 \hline
1,6,7,8 & 1,3 \\ \hline
2,3,4,5 & 2,4 \\ \hline
\end{array}} \ar[dd]^{Q_{\ll}}\\
& & & \\
{\begin{array}{|c|c|}
 \hline
1,2 & 1 \\ \hline
1,1 & 3 \\ \hline
1,1 & 1 \\ \hline
2,2 & 2 \\ \hline
\end{array}} \ar[rrr]^(.3){\pi_{\mu}(V)} & & & {\begin{array}{|c|c|}
 \hline
1,1,1,2 & 1,3 \\ \hline
1,1,2,2 & 1,2 \\ \hline
\end{array}
+
\begin{array}{|c|c|}
 \hline
1,1,1,2 & 1,1 \\ \hline
1,1,2,2 & 2,3 \\ \hline
\end{array}
+
\begin{array}{|c|c|}
 \hline
1,2,2,2 & 1,2 \\ \hline
1,1,1,1 & 1,3 \\ \hline
\end{array}}\\
}
\]
\end{example}

Restricting \ref{eq:pimuspecialization} to the $\ll$-highest weight spaces, we obtain a commutative diagram
\[
 \xymatrix{
 {\rm{hwt}_{\ll}(U_r^{\d})} \ar@<1ex>[rrr]^{Q_{\ll}} \ar[d]_{\pi_{\mu}} & & & {\rm{hwt}_{\ll}(U_r^{\d}(V)) \ar@<1ex>[lll]^{P_{\ll}} \ar[d]^{\pi_{\mu}(V)}}\\
 {\rm{hwt}_{\ll}(U_{\mu}^{\d})} \ar@<1ex>[rrr]^{Q_{\ll}} & & & {\rm{hwt}_{\ll}(U_{\mu}^{\d}(V)) \ar@<1ex>[lll]^{P_{\ll}}}\\
}
\]
where all the horizontal maps are isomorphisms. This shows that the $\ll$-highest weight spaces of the kernels of $\pi_{\mu}$ and $\pi_{\mu}(V)$ get identified via the polarization and specialization maps, and therefore the same is true for $I_r^{\d}$ and $I_r^{\d}(V)$: the generic multi-prolongations and multi-prolongations correspond to each other via polarization and specialization. We summarize the conclusions of this section in the following

\begin{proposition}\label{prop:genericspecial}
 The polarization and specialization maps $P_{\ll}$ and $Q_{\ll}$ restrict to maps between generic flattenings and flattenings, inducing inverse isomorphisms
\[\rm{hwt}_{\ll}(F_{A,B}^{k,r})\simeq\rm{hwt}_{\ll}(F_{A,B}^{k,r}(V)).\]
They also restrict to maps between the kernels of the generic $\pi_{\mu}$'s and the nongeneric ones, inducing inverse isomorphisms
\[\rm{hwt}_{\ll}(\rm{ker}(\pi_{\mu}))\simeq\rm{hwt}_{\ll}(\rm{ker}(\pi_{\mu}(V))).\]
As a consequence, $P_{\ll}$ and $Q_{\ll}$ yield inverse isomorphisms between the $\ll$-highest weight spaces of generic and nongeneric multi-prolongations
\[\rm{hwt}_{\ll}(I_r^{\d})\simeq\rm{hwt}_{\ll}(I_r^{\d}(V)).\]
\end{proposition}

It follows that in order to show that flattenings coincide with multi-prolongations for the variety of secant lines to a Segre-Veronese variety (Theorem \ref{thm:main}), it suffices to prove their equality in the generic setting.

\section{The secant line variety of a Segre-Veronese variety}\label{chap:secantline}

This section is based on the techniques developed in the preceding one. We use the reduction to the ``generic'' situation to work out the analysis of the equations and coordinate rings of secant varieties of Segre-Veronese varieties in the first new interesting case, that of varieties of secant lines. We show how in the case of the secant line variety $\s_2(X)$ of a Segre-Veronese variety $X$, the combinatorics of tableaux can be used to show that the ``generic equations'' coincide with the $3\times 3$ minors of ``generic flattenings''. In particular, we confirm a conjecture of Garcia, Stillman and Sturmfels, which constitutes the special case when $X$ is a Segre variety. We also obtain the representation theoretic description of the homogeneous coordinate ring of $\s_2(X)$, which in particular can be used to compute the Hilbert function of $\s_2(X)$. In the special cases when $\s_2(X)$ coincides with the ambient space, we obtain the decomposition into irreducible representations of certain plethystic compositions. Section \ref{sec:results} describes the statements of our results, while Section \ref{sec:proofs} contains the details of the proofs.

\subsection{Main result and consequences}\label{sec:results}

The main result of our thesis is the description of the generators of the ideal of the variety of secant lines to a Segre-Veronese variety, together with the decomposition of its coordinate ring as a sum of irreducible representations.

\begin{theorem}\label{thm:main} Let $X=SV_{d_1,\cdots,d_n}(\bb{P}V_1^*\times\bb{P}V_2^*\times\cdots\times\bb{P}V_n^*)$ be a Segre-Veronese variety, where each $V_i$ is a vector space of dimension at least $2$ over a field $K$ of characteristic zero. The ideal of $\s_2(X)$ is generated by $3\times 3$ minors of flattenings, and moreover, for every nonnegative integer $r$ we have the decomposition of the degree $r$ part of its homogeneous coordinate ring
\[K[\s_2(X)]_r=\bigoplus_{\substack{\ll=(\ll^1,\cdots,\ll^n)\\ \ll^i\vdash rd_i}}(S_{\ll^1}V_1\o\cdots\o S_{\ll^n}V_n)^{m_{\ll}},\]
where $m_{\ll}$ is obtained as follows. Set 
\[f_{\ll}=\max_{i=1,\cdots,n}\left\lceil\frac{\ll_2^i}{d_i}\right\rceil,\quad e_{\ll}=\ll^1_2+\cdots+\ll^n_2.\]
If some partition $\ll^i$ has more than two parts, or if $e_{\ll}<2f_{\ll}$, then $m_{\ll}=0$. If $e_{\ll}\geq r-1$, then $m_{\ll}=\lfloor r/2\rfloor-f_{\ll}+1$, unless $e_{\ll}$ is odd and $r$ is even, in which case $m_{\ll}=\lfloor r/2\rfloor-f_{\ll}$. If $e_{\ll}<r-1$ and $e_{\ll}\geq 2f_{\ll}$, then $m_{\ll}=\lfloor (e_{\ll}+1)/2\rfloor-f_{\ll}+1$, unless $e_{\ll}$ is odd, in which case  $m_{\ll}=\lfloor (e_{\ll}+1)/2\rfloor-f_{\ll}$. 
\end{theorem}

As a consequence, we derive the conjecture by Garcia, Stillman and Sturmfels, concerning the equations of the secant line variety of a Segre variety.

\begin{corollary}\label{cor:gss}
 The GSS conjecture (Conjecture \ref{conj:gss}) holds, namely the ideal of the variety of secant lines to a Segre product of projective spaces is generated by $3\times 3$ minors of flattenings. 
\end{corollary}

\begin{proof}
 This is the special case of the first part of Theorem \ref{thm:main} when $d_1=d_2=\cdots=d_n=1$.
\end{proof}

Combining Theorem \ref{thm:main} with known dimension calculations for secant varieties of Segre and Veronese varieties, we obtain two interesting plethystic formulas. We do not claim that these formulas are new: since all the vector spaces involved have dimension two, the representation theory of $\sl_2$ can be also used to deduce them. However, we hope that the simple idea we present, together with a generalization of the last part of Theorem \ref{thm:main} to higher secant varieties, would yield new plethystic formulas for decomposing Schur functors applied to tensor products of representations.

\begin{corollary}\label{cor:pletysm3factors}
a) Let $V_1,V_2,V_3$ be vector spaces of dimension two over a field $K$ of characteristic zero, and let $r$ be a positive integer. We have the decomposition
\[\Sym^r(V_1\o V_2\o V_3)=\bigoplus_{\substack{\ll=(\ll^1,\ll^2,\ll^3)\\ \ll^i\vdash r}}(S_{\ll^1}V_1\o S_{\ll^2}V_2\o S_{\ll^3}V_3)^{m_{\ll}},\]
where $m_{\ll}$ is obtained as follows. Set 
\[f_{\ll}=\max\{\ll_2^1,\ll_2^2,\ll_2^3\},\quad e_{\ll}=\ll^1_2+\ll^2_2+\ll^3_2.\]
If some partition $\ll^i$ has more than two parts, or if $e_{\ll}<2f_{\ll}$, then $m_{\ll}=0$. If $e_{\ll}\geq r-1$, then $m_{\ll}=\lfloor r/2\rfloor-f_{\ll}+1$, unless $e_{\ll}$ is odd and $r$ is even, in which case $m_{\ll}=\lfloor r/2\rfloor-f_{\ll}$. If $e_{\ll}<r-1$ and $e_{\ll}\geq 2f_{\ll}$, then $m_{\ll}=\lfloor (e_{\ll}+1)/2\rfloor-f_{\ll}+1$, unless $e_{\ll}$ is odd, in which case  $m_{\ll}=\lfloor (e_{\ll}+1)/2\rfloor-f_{\ll}$. 

b) Let $V_1,V_2$ be vector spaces of dimension two over a field $K$ of characteristic zero, let $r$ be a positive integer and let $\mu=(\mu_1,\mu_2)$ be a partition of $r$ with at most two parts. We have the decomposition
\[S_{\mu}(V_1\o V_2)=\bigoplus_{\substack{\ll=(\ll^1,\ll^2)\\ \ll^i\vdash r}}(S_{\ll^1}V_1\o S_{\ll^2}V_2)^{m_{\ll}},\]
with $m_{\ll}=m_{(\ll^1,\ll^2,\mu)}$, where $m_{(\ll^1,\ll^2,\mu)}$ is as defined in part a).
\end{corollary}

\begin{proof}
 Part a) follows from the fact that the secant line variety of a $3$-factor Segre variety $X$ has the expected dimension, namely $2\cdot\rm{dim}(X)+1$. In the case we are interested in $X=\rm{Seg}(\bb{P}^1\times\bb{P}^1\times\bb{P}^1)$ has dimension $3$ and is a subvariety of $\bb{P}^{2\cdot 2\cdot 2-1}=\bb{P}^7$, so $\s_2(X)$ fills in the whole space. This means that the coordinate ring of $\s_2(X)$ and $\bb{P}^7$ coincide, i.e.
\[K[\s_2(X)]=\Sym(V_1\o V_2\o V_3),\]
and therefore we can use the description of Theorem \ref{thm:main} to compute
\[K[\s_2(X)]_r=\Sym^r(V_1\o V_2\o V_3).\]

 As for part b), let $V_3$ be another vector space of dimension two. Part a) tells us how to decompose
\[\Sym^r(V_1\o V_2\o V_3)\]
in general. On the other hand, regarding $V_1\o V_2\o V_3$ as the tensor product between the vector spaces $V_1\o V_2$ and $V_3$, we can use Cauchy's formula to obtain
\[\Sym^r(V_1\o V_2\o V_3)=\bigoplus_{\mu\vdash r}S_{\mu}(V_1\o V_2)\o S_{\mu} V_3.\]
Now the desired formula for the multiplicity of the irreducible representations occurring in $S_{\mu}(V_1\o V_2)$ follows by combining the formula from part a) with the Cauchy formula depicted above. 
\end{proof}

\begin{corollary}\label{cor:symmetricplethysm}
 Let $V$ be a vector space of dimension two over a field $K$ of characteristic zero. We have the decomposition
\[\Sym^r(\Sym^3(V))=\bigoplus_{\ll\vdash 3r} (S_{\ll}V)^{m_{\ll}},\]
where $m_{\ll}$ is obtained as follows. Set 
\[f_{\ll}=\left\lceil\frac{\ll_2}{3}\right\rceil,\quad e_{\ll}=\ll_2.\]
If $\ll$ has more than two parts, or if $e_{\ll}<2f_{\ll}$ (i.e. $\ll_2=1$), then $m_{\ll}=0$. If $e_{\ll}\geq r-1$, then $m_{\ll}=\lfloor r/2\rfloor-f_{\ll}+1$, unless $e_{\ll}$ is odd and $r$ is even, in which case $m_{\ll}=\lfloor r/2\rfloor-f_{\ll}$. If $e_{\ll}<r-1$ and $e_{\ll}\geq 2f_{\ll}$, then $m_{\ll}=\lfloor (e_{\ll}+1)/2\rfloor-f_{\ll}+1$, unless $e_{\ll}$ is odd, in which case  $m_{\ll}=\lfloor (e_{\ll}+1)/2\rfloor-f_{\ll}$. 
\end{corollary}

\begin{proof}
 This follows from the fact that $\s_2(\rm{Ver}_3(\bb{P}^1))$, the secant line variety of the twisted cubic fills in the space, hence its coordinate ring is
\[\Sym(\Sym^3(V)).\]
Using the description in Theorem \ref{thm:main} with $n=1$, $d_1=3$ and $V=V_1$ of dimension $2$, we obtain the desired formula.
\end{proof}

\subsection{Proof of the main result}\label{sec:proofs}

\begin{proof} We start by outlining the main steps in the proof of Theorem \ref{thm:main}. We fix a sequence of positive integers $\d=(d_1,\cdots,d_n)$ and a positive degree $r$, and let $\r=(rd_1,\cdots,rd_n)$. By Proposition \ref{prop:genericspecial}, it suffices to prove the generic version of the theorem. More precisely, we let
\[F=\sum_{A+B=\d}F_{A,B}^{3,r}\subset U_r^{\d}\]
be the set of generic flattenings, and let $F_i$ denote those generic flattenings with $|A|=i$. (As the rest of the proof will imply, we have $F=F_1+F_2+F_3$; see \cite{rai3x3} or \cite[Chapter 6]{raicu} for more precise results in this direction in the case $n=1$ of the Veronese variety.)

Recall that $I=I_r^{\d}$ denotes the space of generic multi-prolongations of degree $r$ (Definition \ref{def:genericprolongations}), i.e. $I$ is the kernel of the map
\[\pi=\bigoplus_{\mu=(\mu_1,\mu_2)\vdash r}\pi_{\mu}:U=U_r^{\d}\lra\bigoplus_{\mu=(\mu_1,\mu_2)\vdash r}U_{\mu}^{\d}.\]
We have $F\subset I$, by combining Lemma \ref{lem:eqnflat} with Proposition \ref{prop:genericspecial}. We will show that $F=I$ and that the image of $\pi$ decomposes into irreducible $S_{\r}$-representations as
\[\pi(U)=\bigoplus_{\ll\vdash^n\r}[\ll]^{m_{\ll}},\]
where $m_{\ll}$ is as defined in the statement of the theorem.

We list the main steps below. The details will occupy the rest of the section.

\defi{Step 0:} If $\ll$ is an $n$-partition with some $\ll^i$ having at least three parts, then $\rm{hwt}_{\ll}(U)=\rm{hwt}_{\ll}(F)$ (Proposition \ref{prop:1flattenings}), hence $\rm{hwt}_{\ll}(F)=\rm{hwt}_{\ll}(I)$, because $F\subset I\subset U$. Moreover, this also shows that $m_{\ll}=0$.

\defi{Step 1:} We fix an $n$-partition $\ll$ of $\r$ with each $\ll^i$ having at most two parts. We identify each tableau $T$ with a certain graph $G$. We show that graphs containing odd cycles are contained in $F$.

\defi{Step 2:} We show that the $\ll$-highest weight space of $U/F$ is spanned by bipartite graphs that are as connected as possible, i.e. that are either connected, or a union of a tree and some isolated nodes.

\defi{Step 3:} We introduce the notion of \defi{type} associated to a graph $G$ as in \defi{Step 2}, encoding the sizes of the sets in the bipartition of the maximal component of $G$. We show that if $G_1,G_2$ have the same type, then $G_1=\pm G_2\rm{\ (modulo\ }F)$.

\defi{Step 4:} If we let
\[\pi=\bigoplus_{\substack{\mu\vdash r\\ \mu=(a\geq b)}}\pi_{\mu}:U\lra\bigoplus_{\substack{\mu\vdash r\\ \mu=(a\geq b)}}U_{\mu}^n,\]
and if $G_i$ are graphs of distinct types (not contained in $F$), then the elements $\pi(G_i)$ are linearly independent. This suffices to prove that $F$ and the kernel of $\pi$ are the same, i.e. that $F=I$. The formulas for the multiplicities $m_{\ll}$ follow from counting the number of $G_i$'s, i.e. the number of possible types.

\subsubsection{Step 1} We fix an $n$-partition $\ll$ of $\r$ with $\ll^i=(\ll^i_1\geq\ll^i_2\geq 0)$, for $i=1,\cdots,n$. For each $n$-tableau $T$ of shape $\ll$ we construct a graph $G$ with $r$ vertices labeled by the elements of the alphabet $\mc{A}=\{1,\cdots,r\}$ as follows. For each tableau $T^i$ of $T$ and column $\Yvcentermath1 \young(x,y)$ of $T^i$ of length $2$, $G$ has an oriented edge $(x,y)$ which we label by the index $i$. We will often refer to the labels of the edges of $G$ as colors. Note that we allow $G$ to have multiple edges between two vertices (some call such $G$ a multigraph), but at any given vertex there can be at most $d_i$ incident edges of color $i$. Since we think of two $n$-tableaux as being the same if they differ by a permutation of $\mc{A}$, we shall also identify two graphs if they differ by a relabeling of their nodes. Note that a graph $G$ determines an element in $\rm{hwt}_{\ll}(U)$, by considering a tableau $T$ with columns $\Yvcentermath1 \young(x,y)$ for each edge $(x,y)$ of $G$. The order of the columns of $T$ is not determined by $G$, but part (2) of Lemma \ref{lem:relstableaux} states that any such $T$ yields the same element of $\rm{hwt}_{\ll}(U)$. The orientation of the edges of our graphs will be mostly irrelevant: reversing the orientation of an edge of $G=T$ will correspond to changing $G$ to $-G$ (see part (1) of Lemma \ref{lem:relstableaux}). When we talk about connectedness and cycles, we don't take into account the orientation of the edges.

\begin{example} The graph
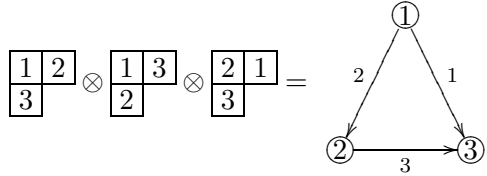
\[\xymatrix@=4pt
{
& & & \ccircle{1} \ar[lldd]_2 \ar[rrdd]^1 & &\\
{\Yvcentermath1 \young(12,3)\o\young(13,2)\o\young(21,3) =} & & & & & \\
& \ccircle{2} \ar[rrrr]_3 & & & & \ccircle{3}\\
}
\]
is connected and has a cycle of length $3$, while
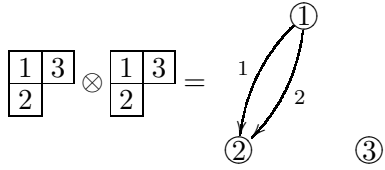
\[\xymatrix@=4pt
{
& & & \ccircle{1} \ar@/_/[lldd]_1 \ar@/^/[lldd]^2 & &\\
{\Yvcentermath1 \young(13,2)\o\young(13,2) =} & & & & & \\
& \ccircle{2} & & & & \ccircle{3}\\
}
\]
is disconnected and has a cycle of length $2$.
\end{example}

From now on we work modulo $F$, and more precisely, inside the $\ll$-highest weight space of $(U/F)$. This space is generated by the graphs described above. The main result of \defi{Step 1} is

\begin{proposition}\label{propoddcyc} If $G$ has an odd cycle, then $G=0$ (i.e. $G$ is in $F$).
\end{proposition}

We first need to establish some fundamental relations, that will be used throughout the rest of the proof.

\begin{lemma}\label{lem:fundrels} The following relations between tableaux/graphs hold (see the interpretation below)

a) $\Yvcentermath1\young(x,y)=-\young(y,x)$, in particular $\Yvcentermath1\young(x,x)=0$.

b) $\Yvcentermath1\young(xz,y)=\young(xy,z)+\young(zx,y)$.

c) $\Yvcentermath1\young(xz,y)\o\young(xy,z)=\young(xy,z)\o\young(xz,y)$.

d) $\Yvcentermath1\young(xz,y)\o\young(xz,y)\o\young(xy,z)=\young(xz,y)\o\young(xy,z)\o\young(xy,z)$.
\end{lemma}

\emph{Interpretation:} For an expression $E=\sum_T a_T\cdot T$, where the $T$'s are $n$-tableaux of shape $\ll$, we say that $E=0$ if
\[\sum_T a_T\cdot T \in \ F\subset U.\]
If all the $n$ tableaux occurring in the expression $E$ contain the same $n$-subtableau $S$, then we suppress $S$ entirely from the notation (see also the comment in part (3) of Lemma \ref{lem:relstableaux}).

\begin{example} One interpretation of part b) of Lemma \ref{lem:fundrels} could be that
\[\Yvcentermath1\young(ab,cd)\o\young(xzt,y)=\young(ab,cd)\o\young(xyt,z)+\young(ab,cd)\o\young(zxt,y),\]
for any $\{a,b,c,d\}=\{x,y,z,t\}=\{1,2,3,4\}$. The $2$-subtableau $S$ is in this case
\[{S=\Yvcentermath1\young(ab,cd)\ \o\ }\young(t)\ .\]
\end{example}

\begin{proof}[Proof of Lemma \ref{lem:fundrels}] a) is part (1) of Lemma \ref{lem:relstableaux}.

b) follows from part (3) of the same lemma (since all columns of our tableaux have size at most two).

c) We have
\[\Yvcentermath1 \young(\ccx\ccy,\ccz)\o\young(\ccx\ccz,\ccy)=\sum_{\s\in S_3}\rm{sgn}(\s)\cdot\s\left(\young(xy,z)\o\young(xz,y)\right)=0,\]
(because the left hand side is contained in $F_2$). Using parts a) and b) repeatedly, we can express everything in terms of $\Yvcentermath1\young(xy,z)$ and $\Yvcentermath1\young(xz,y)$, and after simplifications, the above equation becomes
\[\Yvcentermath1 3\cdot\left(\young(xy,z)\o\young(xz,y)-\young(xz,y)\o\young(xy,z)\right)=0.\]

d) Part c) states that any tensor expression in $a=\Yvcentermath1\young(xz,y)$ and $b=\Yvcentermath1\young(xy,z)$ does not depend on the order in which $a$ and $b$ appear, so we can think of the pure tensors in $a,b$ as commuting monomials in $a,b$. Writing $\Yvcentermath1\young(yx,z)=b-a$, we can translate
\[\Yvcentermath1 \young(\ccx\ccz,\ccy)\o\young(\ccx\ccz,\ccy)\o\young(\ccx\ccy,\ccz)=\sum_{\s\in S_3}\rm{sgn}(\s)\cdot\s\left(\young(xz,y)\o\young(xz,y)\o\young(xy,z)\right)=0,\]
into
\[a^2b-a^2(b-a)+(a-b)^2b-b^2a-(b-a)^2a+b^2(a-b)=0,\]
which simplifies to $3(a^2b-ab^2)=0$, i.e. $a^2b=ab^2$, or
\[\Yvcentermath1\young(xz,y)\o\young(xz,y)\o\young(xy,z)=\young(xz,y)\o\young(xy,z)\o\young(xy,z).\]
\end{proof}

\begin{corollary}\label{corodd12} If $G$ is a graph having a connected component $H$ consisting of two nodes joined by an odd number of edges, then $G=0$.
\end{corollary}

\begin{proof} Interchanging the labels of the two nodes of $H$ preserves $G$, but by part a) of Lemma \ref{lem:fundrels}, it also transforms $G$ into $(-1)^e G$, where $e$ is the number of edges in $H$. Since $e$ is odd, $G=0$.
\end{proof}

\begin{corollary}\label{corodd13} If $G$ is a graph containing cycles of length $1$ or $3$, then $G=0$.
\end{corollary}

\begin{proof} If $G$ has a cycle of length $1$, this follows from part a) of Lemma \ref{lem:fundrels}. If $G$ has a cycle of length $3$, we may assume this cycle is $C=\ccone\to\cctwo\to\ccthree\to\ccone$. We have several cases to analyze, depending on the colors of the edges in this cycle.

If the edges in $C$ have distinct colors, we need to prove that
\[\Yvcentermath1\young(13,2)\o\young(21,3)\o\young(12,3)=0.\]
We have by part b) of Lemma \ref{lem:fundrels} applied to the middle tableau that
\[\Yvcentermath1\young(13,2)\o\young(21,3)\o\young(12,3)=\young(13,2)\o\young(12,3)\o\young(12,3)-\young(13,2)\o\young(13,2)\o\young(12,3)=0,\]
where the last equality is part d) of the same lemma.

If the edges of $C$ have the same color, we need to prove that
\[\Yvcentermath1\young(112,233)=0.\]
We have
\[\Yvcentermath1 0=\young(\ccone\ccone\cctwo,\cctwo\ccthree\ccthree)=\young(112,233)-\young(221,133)-\young(332,211)-\young(113,322)+\young(223,311)+\young(331,122)\]
\[\Yvcentermath1 =\young(112,233)+\young(121,233)+\young(211,332)+\young(112,323)+\young(211,323)+\young(121,332)=6\cdot\young(112,233),\]
where the penultimate equality follows from skew-symmetry on rows, while the last one follows from part (2) of Lemma \ref{lem:relstableaux}.

Finally, suppose that the edges of $C$ have two colors, say $(1,2)$ and $(1,3)$ have the same color. We need to prove that
\[\Yvcentermath1\young(1123,23)\o\young(21,3)=0.\]
As in the preceding case,
\[\Yvcentermath1 0=\young(\ccone\ccone\cctwo\ccthree,\cctwo\ccthree)\o\young(\cctwo\ccone,\ccthree)=\young(1123,23)\o\young(21,3)-\young(2213,13)\o\young(12,3)-\young(3321,21)\o\young(23,1)\]
\[\Yvcentermath1 -\young(1132,32)\o\young(31,2)+\young(2231,31)\o\young(32,1)+\young(3312,12)\o\young(13,2)=6\cdot\young(1123,23)\o\young(21,3),\]
where the last equality follows by utilizing repeatedly parts a) and c) of \ref{lem:fundrels}. For example, we have for the second term that
\[\Yvcentermath1\young(2213,13)\o\young(12,3)=-\young(1213,23)\o\young(12,3)=-\young(1123,23)\o\young(21,3),\]
where the last equality follows by applying part c) of \ref{lem:fundrels} in the form
\[\Yvcentermath1\young(yz,x)\o\young(zy,x)=\young(zy,x)\o\young(yz,x),\]
with
\[\Yvcentermath1\young(yz,x)=\young(21,3),\quad \young(zy,x)=\young(12,3).\]
\end{proof}

\begin{corollary}\label{corswitch23} If an $n$-tableau $T$ contains the columns $C_1=\Yvcentermath1\young(x,y)$ and $\Yvcentermath1 C_2=\young(x,z)$, and $T'$ is obtained from $T$ by interchanging two boxes $\boxed{y}$ and $\boxed{z}$ from the same tableau $T^i$ of $T$, and not contained in any of $C_1,C_2$, then $T=T'$ (modulo $F$).
\end{corollary}

\begin{proof} If $\Yvcentermath1\young(y,z)$ is a column of $T^i$ then $T$ contains a triangle, hence $T=0$. Since interchanging $y$ and $z$ transforms $T$ into $T'=-T=0$, it follows that $T=T'$. We can assume then that $y$ and $z$ don't lie in the same column of $T^i$. If they both belong to columns of size one of $T^i$, then interchanging them preserves $T$ (see part (2) of Lemma \ref{lem:relstableaux}). Otherwise we may assume that $y$ belongs to a column of size two in $T^i$, hence we have the relation
\[\Yvcentermath1\young(yz,*)=\young(y*,z)+\young(zy,*)=\young(zy,*),\]
where the last equality follows from the fact that any tableau containing $C_1$, $C_2$ and $\Yvcentermath1\young(y,z)$ is a graph containing a triangle, i.e. it is zero (Corollary \ref{corodd13}).
\end{proof}

\begin{proof}[Proof of Proposition \ref{propoddcyc}] We show that a graph $G$ (with corresponding tableau $T$) containing an odd cycle of length at least $5$ is a linear combination of graphs with shorter odd cycles. The conclusion then follows by induction from Corollary \ref{corodd13}. Suppose that $C:\ccirc{1}\to \ccirc{2}\to\cdots\to \ccirc{k}\to \ccirc{1}$ is an odd cycle in $G$, with $k\geq 5$. We denote by $E_i$ the edge $(i,i+1)$ ($E_k=(k,1)$).

Let's assume first that there are two consecutive edges of $C$ of the same color, say $E_1$ and $E_2$ have color $1$. If not all edges of $C$ have color $1$, we may assume that $E_3$ has color $2$, so that $T$ contains the subtableau
\[\Yvcentermath1\young(12,23)\o\young(3,4).\]
Since $E_1,E_2$ have color $1$, it follows that $d_1\geq 2$, hence there are at least two $4$'s in $T^1$. One of them is thus not contained in $E_5$, and therefore in none of the edges of $C$. We apply Corollary \ref{corswitch23} with $C_1,C_2$ the columns corresponding to $E_2,E_3$, $y=2$ and $z=4$. We can thus interchange the $2$ in $E_1$ with a $4\in T^1$ not in any $E_i$, obtaining an $n$-tableau $T'=T$, with $T'$ containing the cycle $\ccirc{1}\to\ccirc{4}\to\ccirc{5}\to\cdots\to\ccirc{k}\to\ccirc{1}$ of length $k-2$.

If all the $E_i$'s have color $1$, $T$ contains the subtableau
\[\Yvcentermath1 S=\young(1234\cdots k,2345\cdots 1)\subset T^1.\]
If there is an edge $(3,4)$ of $G$ with color different from $1$, then we can replace $E_3$ by that edge and apply the previous case. If $d_1>2$ then $T^1$ has a $4$ not contained in any $E_i$, so we can again use the argument from the previous paragraph. Suppose now that $d_1=2$. The proof of Corollary \ref{corswitch23} shows that we can interchange $3$ with $4$ in all $T^i$'s ($i\neq 1$), modulo tableaux containing $S$ and an edge $(3,4)$ of color different from $1$. But these we know are zero (modulo $F$) by the argument above, so we can write $T=T'$ where $T'$ is obtained from $T$ by interchanging all $3$'s and $4$'s in $T^i$ for $i\geq 2$. We now use the relation
\[\Yvcentermath1 \young(1234,2345)=\young(1234,2543)+\young(1233,2445),\]
to write
\[T'=T''+T''',\]
where $T''$ contains the cycle $\ccirc{1}\to\ccirc{2}\to\ccirc{5}\to\cdots\to\ccirc{k}\to\ccirc{1}$ of length $k-2$, and $T'''$ differs from $T$ by interchanging all the $3$'s and $4$'s in $T$, and doing a column transposition in the column of $E_3$. This shows that $T=T'=0-T$, hence $T=0$.

Finally, we assume that no two consecutive edges have the same color. Since the cycle is odd, we can find three consecutive edges with distinct colors, say $E_1$, $E_2$ and $E_3$, with colors $1,2$ and $3$ respectively. By Corollary \ref{corswitch23}, we have
\[\Yvcentermath1 T = \young(14,2)\o\young(2,3)\o\young(3,4)=\young(12,4)\o\young(2,3)\o\young(3,4).\]
If the edge $E_4$ in $C$ doesn't have color $1$, then it survives after interchanging $2$ and $4$ as above, hence $T$ is equal with a graph containing the odd cycle $\ccirc{1}\to \ccirc{4}\to \ccirc{5}\to\cdots\to \ccirc{k}\to \ccirc{1}$ of length $k-2$.

Suppose now that $E_4$ has color $1$. If the edge $E_5$ doesn't have color $2$, then we may repeat the above argument replacing the edges $E_1$, $E_2$ and $E_3$ with $E_2$, $E_3$ and $E_4$ respectively. Otherwise, $T$ contains the subtableau (with $*=6$ if $k>5$ and $*=1$ if $k=5$)
\[\Yvcentermath1\young(14,25)\o\young(25,3*)\o\young(35,4)=\young(12,45)\o\young(25,3*)\o\young(35,4)=\young(12,45)\o\young(25,3*)\o\young(53,4),\]
where the first equality follows by interchanging $2$ and $4$ in the first factor, while the last one follows by interchanging $3$ and $5$ in the last factor (in both cases we apply Corollary \ref{corswitch23}). It follows that $T$ is equal to a graph containing the odd cycle $\ccirc{1}\to \ccirc{4}\to \ccirc{5}\to\cdots\to \ccirc{k}\to \ccirc{1}$ of length $k-2$, concluding the proof.
\end{proof}

\subsubsection{Step 2} We first translate the relations in part b) of Lemma \ref{lem:fundrels} into \defi{basic operations} on graphs. We start with the following

\begin{definition} A node $\ccirc{j}$ is said to be $i$-saturated if there are $d_i$ edges of color $i$ incident to $\ccirc{j}$.
\end{definition}

\begin{remark}[Basic operations]\label{rembasicop} Let $G$ be a graph containing an edge $(1,2)$ of color $1$. The following relations hold:
\begin{enumerate}

\item If the vertex $\ccirc{3}$ is not $1$-saturated, then $\Yvcentermath1\young(13,2)=\young(12,3)+\young(31,2)$ becomes
\[\xymatrix@=4pt
{
 & \ccircle{1} \ar[ldd]_1 & & & & \ccircle{1} \ar[rdd]^1 & & & & \ccircle{1} & \\
 & & &  = & & & & + & & & \\
\ccircle{2} & & \ccircle{3} & & \ccircle{2} & & \ccircle{3} & & \ccircle{2} & & \ccircle{3} \ar[ll]^1\\
}
\]

\item If $G$ has an edge $(3,4)$ of color $1$, then $\Yvcentermath1\young(13,24)=\young(12,34)+\young(13,42)$ becomes
\[\xymatrix@=4pt
{
\ccircle{1} \ar[dd]_1 & \ccircle{3} \ar[dd]^1 & & \ccircle{1} \ar[rr]^1 & & \ccircle{3} & & \ccircle{1} \ar[rrdd]_<<<1 & & \ccircle{3} \ar[lldd]^<<<1 \\
 & & = & & & & + & & \\
\ccircle{2} & \ccircle{4} & & \ccircle{2} \ar[rr]_1 & & \ccircle{4} & & \ccircle{2} & & \ccircle{4} \\
}
\]
\end{enumerate}
\end{remark}

\begin{proposition}\label{prop:MCBgenerate} Let $\ll$ be as before, and let
\[e_{\ll}=\sum_{i=1}^n\ll^i_2.\]
If $e_{\ll}\geq r-1$, then $\rm{hwt}_{\ll}(U/F)$ is spanned by connected graphs. If $e_{\ll}<r-1$, then $\rm{hwt}_{\ll}(U/F)$ is spanned by graphs $G$ that consist of a tree, together with a collection of isolated nodes.
\end{proposition}

\begin{proof} We first show that if $G$ has two connected components $H_1,H_2$ with $H_1$ containing a cycle, then we can write $G=G_1+G_2$, where $G_1$ and $G_2$ are graphs obtained from $G$ by joining the components $H_1,H_2$ together.

Consider an edge $(1,2)$ contained in a cycle of $H_1$, having say color $1$. Consider a node $\ccirc{3}$ of $H_2$ and suppose first it is not $1$-saturated. Using the first basic operation of Remark \ref{rembasicop}, we get that $G=G_1+G_2$, where $G_1,G_2$ are obtained from $G$ by connecting $H_2$ to $H_1$ via an edge of color $1$. If $\ccirc{3}$ is $1$-saturated, then in particular there exists at least one edge, say $(3,4)$, of color $1$ in $H_2$. The second basic operation of Remark \ref{rembasicop} yields $G=G_1+G_2$, where $G_1,G_2$ are obtained from $G$ by connecting $H_1$ and $H_2$ via two edges of color $1$.

If $e_{\ll}\geq r-1$, then $G$ will contain cycles as long as it is not connected, so iterating the above procedure, we can write $G$ as a linear combination of connected graphs.

If $e_{\ll}<r-1$, then the above argument reduces the problem to the case when $G$ is a union of trees, some of which may be isolated nodes. We show that if $G$ has at least two components that are not nodes, then $G=G_1+G_2$, where $G_1,G_2$ are unions of trees, and the sizes of the largest components of $G_1,G_2$ are strictly larger than the size of the largest component of $G$. Induction on the size of the largest component of $G$ concludes then the proof of the proposition.

Let $H_1$ be the largest component of $G$, and let $H_2$ be another component which isn't a node. If $H_2$ has only one edge, then $G=0$ by Corollary \ref{corodd12}. Consider a leaf of $H_1$, say $\ccirc{3}$, and assume first that all edges in $H_2$ have the same color, say $1$. Since $H_2$ has more than one edge and is connected, it must have a vertex with at least two incident edges of color $1$, i.e. $d_1\geq 2$. This means that $\ccirc{3}$ is not $1$-saturated. Let $(1,2)$ be an edge of $H_2$ (of color $1$). The first basic operation of Remark \ref{rembasicop} shows that $G=G_1+G_2$, where $G_1,G_2$ are obtained from $G$ by expanding its largest component.

Assume now that the edges in $H_2$ have at least two colors, and that the edge incident to $\ccirc{3}$ has color $2$. Let $(1,2)$ be an edge of $H_2$ of color different from $2$, say $1$. $\ccirc{3}$ is not $1$-saturated, thus we can use the first basic operation of Remark \ref{rembasicop} as in the preceding case.
\end{proof}

\subsubsection{Step 3} Combining \defi{Step 1} with \defi{Step 2} we get that, depending on the $n$-partition $\ll$, $\rm{hwt}_{\ll}(U/F)$ is spanned either by connected graphs without odd cycles, or by graphs consisting of a tree and some isolated nodes. We call these graphs \defi{maximally connected bipartite} (MCB) graphs. For an MCB-graph $G$, the maximal connected component admits an essentially unique \defi{bipartition} of its vertex set into subsets $A,B$ of sizes $a\geq b$ (i.e. vertices in the same subset $A$ or $B$ are not connected by an edge). We say that $G$ has \defi{type} $(a,b;\ll)$ (or just $(a,b)$ when $\ll$ is understood), and that it is \defi{canonically oriented} if all the edges have source in $A$ and target in $B$ (when $a=b$, there are two canonical orientations). We have the following

\begin{proposition}\label{propMCB} If $G_1,G_2$ are canonically oriented MCB-graphs of type $(a,b)$, then $G_1=G_2$.
\end{proposition}

We first need to refine the relations of Remark \ref{rembasicop}:

\begin{remark}[Refined basic operations]\label{remrefbasicop} Suppose that $G$ is an MCB-graph with vertex bipartition $A\sqcup B$ as above.

\begin{enumerate}

\item Assume that $\ccirc{3}$ is not $1$-saturated, $(1,2)$ is an edge of color $1$, and $\ccirc{1},\ccirc{3}$ belong to $A$. If $\ccirc{1},\ccirc{3}$ are contained in the same connected component of the graph obtained from $G$ by removing the edge $(1,2)$, then
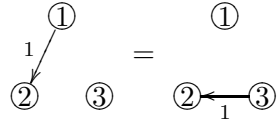
\[\xymatrix@=4pt
{
 & \ccircle{1} \ar[ldd]_1 & & & & \ccircle{1}  &\\
 & & &  = & & &\\
\ccircle{2} & & \ccircle{3} & & \ccircle{2} & & \ccircle{3} \ar[ll]^1\\
}
\]
This follows from the fact that the above conditions guarantee that the term that was left out from the first basic operation of Remark \ref{rembasicop} has an odd cycle, and hence equals $0$ by Proposition \ref{propoddcyc}.

\item Assume that $(1,2)$ and $(3,4)$ are edges of color $1$, $\ccirc{1},\ccirc{3}\in A$ and $\ccirc{2},\ccirc{4}\in B$, and either $\ccirc{1}$ and $\ccirc{3}$, or $\ccirc{2}$ and $\ccirc{4}$ are in the same connected component of the graph obtained from $G$ by removing the edges $(1,2)$ and $(3,4)$. Then
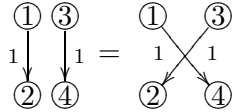
\[\xymatrix@=4pt
{
\ccircle{1} \ar[dd]_1 & \ccircle{3} \ar[dd]^1 &  & \ccircle{1} \ar[rrdd]_<<<1 & & \ccircle{3} \ar[lldd]^<<<1 \\
 & & = & & \\
\ccircle{2} & \ccircle{4} & & \ccircle{2} & & \ccircle{4} \\
}
\]
As above, the missing term from the second basic operation has an odd cycle, and hence equals $0$.
\end{enumerate}
\end{remark}

\begin{proof}[Proof of Proposition \ref{propMCB}] We prove by induction on $e_{\ll}$ (the number of ``edges'' of $\ll$), that it is possible to get from $G_1$ to $G_2$ via a series of refined basic operations. If $e_{\ll}=0$, there is nothing to prove. Suppose now that $e_{\ll}>0$.

We call an edge $E$ of an MCB-graph $G$ \defi{nondisconnecting} if the graph obtained from $G$ by removing $E$ is still an MCB-graph. More explicitly, if $e_{\ll}\geq r$, then $E$ must be contained in a cycle of $G$, and if $e_{\ll}<r$, then one of the endpoints of $E$ must be a leaf of $G$.

We will prove that for any nondisconnecting edge $E_2$ of $G_2$ of color $c$, there exist a sequence of refined basic operations which transforms $G_1$ into a new graph $\hat{G}_1$ having a nondisconnecting edge $E_1$ of color $c$, such that the graphs $G_1'$ and $G_2'$ obtained from $\hat{G}_1$ and $G_2$ by removing the edges $E_1$ and $E_2$ have the same type. Assuming this, by induction we can find a series of refined basic operations that transform $G_1'$ into $G_2'$. We lift this sequence of operations to $\hat{G}_1$ as follows: the refined basic operations of type (2) are performed just as if the edge $E_1$ was not contained in $\hat{G}_1$, as well as the operations of type (1) that don't transform an edge $E'$ of color $c$ into one that's incident to $E_1$; the operations of type (1) involving an edge $E'$ of color $c$ that gets transformed into an edge incident to $E_1$ are replaced by operations of type (2) involving $E'$ and $E_1$. It is clear that $E_1$ remains nondisconnecting along the process, so we end up with the graphs $G_1''$ and $G_2$ that coincide after removing the nondisconnecting edges $E_1$ and $E_2$ of color $c$. At most two more refined operations of type (2) (that correspond to correcting the positions of the endpoints of $E_1$) are then sufficient to transform $G_1''$ into $G_2$, concluding the proof.

We now show that if $e_{\ll}\geq r$ and $G_1$ has an edge $E_1$ of color $c$, then we can find a refined basic operation that makes $E_1$ nondisconnecting. Suppose that $E_1$ is disconnecting, and let $H_1,H_2$ be the connected components of the graph obtained from $G_1$ by removing the edge $E_1$. One of $H_1,H_2$ must contain a cycle, say $H_1$, and let $O$, $Y$ be consecutive edges of this cycle, of colors $o(range)$ and $y(ellow)$ (note that $o$ might coincide with $y$). If $H_2$ has a node $N$ that is not $o$-saturated or not $y$-saturated, then a refined operation of type (1) involving the node $N$ (as $\ccirc{3}$) and one of the edges $O$, $Y$ (as the edge $(1,2)$) will make $E_1$ a nondisconnecting edge. Otherwise, if every vertex of $H_2$ is both $o$- and $y$-saturated, then there exists a cycle in $H_2$ consisting of edges of colors $o$ and $y$ (if $o=y$, then since $O,Y$ are incident edges of color $o$, it means that $d_o\geq 2$, in particular any $o$-saturated node has at least two incident edges of color $o$; if $o\neq y$, then any $o$- and $y$-saturated node has at least one $o$- and one $y$- incident edge; in both cases, the nodes in $H_2$ have at least two incident edges, so we can find a cycle as stated). A refined basic operation of type (2) involving an $o$-edge on this cycle and $O$ (or an $y$-edge and $Y$) will make $E_1$ nondisconnecting.

Finally, if $e_{\ll}<r$ and $G_2$ has a nondisconnecting edge $E_2$ of color $c$, then we have claimed that we can find a sequence of refined basic operations that transforms $G_1$ into a graph $\hat{G}_1$ containing a nondisconnecting edge $E_1$ of color $c$, and moreover $\hat{G}_1-E_1$ and $G_2-E_2$ have the same type. We may assume that $e_{\ll}=r-1$, by removing the isolated nodes of $G_1$ and $G_2$. Suppose that the graphs $G_i$ have vertex bipartitions $A_i\sqcup B_i$, with $|A_i|=a$, $|B_i|=b$, and that $E_2=(x,y)$, with $\ccirc{y}\in B_2$ a leaf of $G_2$. This means that the graph $G_2$, and hence also $G_1$, has at most $(b-1)\cdot d_c+1$ edges of color $c$, and for any color $c'\neq c$, it has at most $(b-1)\cdot d_{c'}$ edges of color $c'$. In particular, for any color $c'\neq c$, there exists a node in $B_1$ which is not $c'$-saturated. Consider an edge $E=(u,v)$ of color $c$ in $G_1$, with $\ccirc{u}\in A_1$, $\ccirc{v}\in B_1$. Let $H_1, H_2$ be the connected components of $G_1-E$ containing $u$ and $v$ respectively. We prove by descending induction on the size of $H_2$ that we can make $E$ nondisconnecting, with its endpoint in $B_1$ being a leaf.

If $H_2=\{\ccirc{v}\}$ then $E$ is nondisconnecting. More generally, if $H_2\cap B_1=\{\ccirc{v}\}$, then we may assume that all the edges in $H_2$ have color $c$. If $E'$ is an edge of $H_2$ of color $c'\neq c$ (see the second transformation in Example \ref{ex:basicopsgraphs} below), then there are at most $(b-1)\cdot d_{c'}-1$ edges of color $c'$ in $H_1$, so that we can find a vertex $C'$ in $H_1$ that is not $c'$-saturated. A refined basic operation of type (1) involving $E'$ and $C'$ decreases the size of $H_2$ by one, so we can conclude by induction. Assume now that the edges in $H_2$ have color $c$. Together with the edge $E$, we get at least two edges of color $c$ outside $H_1$, which means that $H_1$ has at most $(b-1)\cdot d_c-1$ edges of color $c$, i.e. it has a vertex that is not $c$-saturated. We now do a refined basic operation of type (1) as before, involving that vertex and an edge of $H_2$, and conclude by induction.

We may now assume that $|H_2\cap B_1|>1$ (see the first transformation in Example \ref{ex:basicopsgraphs} below). Therefore there exist distinct edges $Y=(u',v)$ of color $y(ellow)$ and $O=(u',v')$ of color $o(range)$ in $H_2$ ($y$ and $o$ might coincide). If $o=c$ then we replace $E$ with $(u',v')$, which decreases the size of $H_2$, so that we can conclude by induction. If there exists a vertex $W\in H_1\cap B_1$ that is not $y$-saturated, then the refined basic operation involving $Y$ and $W$ decreases the size of $H_2$. Likewise, if there exists a vertex $W\in H_1\cap A_1$ that is not $o$-saturated, then the refined basic operation involving $O$ and $W$ also decreases the size of $H_2$. We may therefore assume that all nodes in $B_1\cap H_1$ are $y$-saturated, and those in $A_1\cap H_1$ are $o$-saturated, and show that this leads to a contradiction. If $y=o$, then since $u'$ has two incident edges of color $o$, we must have $d_o\geq 2$. All the nodes of $H_1$ being saturated implies that they have degree at least $d_o\geq 2$, so $H_1$ contains a cycle, which is a contradiction. If $y\neq o$, then $H_1$ must contain at least $|H_1\cap A_1|$ edges of color $o$ (since each vertex in $H_1\cap A_1$ is $o$-saturated) and at least $|H_1\cap B_1|$ edges of color $y$, i.e. $H_1$ contains at least $|H_1|$ edges, hence it can't be a tree.
\end{proof}

\begin{example}\label{ex:basicopsgraphs} Consider the $3$-tableaux 
\[\Yvcentermath1 T_1=\young(3516,24)\o\young(31256,4)\o\young(1534,26),\ T_2=\young(1356,24)\o\young(51234,6)\o\young(3516,24),\]
with corresponding graphs
\[\xymatrix@=15pt{
& \ccircle{1} \ar@{-->}[dr] & \\
& & \ccircle{2} \\
G_1\quad = & \ccircle{3} \ar[ur] \ar@{~>}[dr] & \\
& & \ccircle{4} \\
& \ccircle{5} \ar[ur] \ar@{-->}[dr] & \\
& & \ccircle{6} \\
}
\quad\quad\quad\quad
\xymatrix@=15pt{
& \ccircle{1} \ar[dr] & \\
& & \ccircle{2} \\
{\rm{and}\quad\quad\quad\quad G_2\quad =} & \ccircle{3} \ar@{-->}[ur] \ar[dr] & \\
& & \ccircle{4} \\
& \ccircle{5} \ar@{-->}[ur] \ar@{~>}[dr] & \\
& & \ccircle{6} \\
}
\]
where color $1$ corresponds to $\xymatrix{\ar[r]&}$, color $2$ to $\xymatrix{\ar@{~>}[r]&}$, and color $3$ to $\xymatrix{\ar@{-->}[r]&}$. $G_1$ and $G_2$ are MCB of the same type, and in fact $G_1=0$, since it is the same as the graph obtained by reversing the orientation of its $5$ edges (an odd number), and this equals $-G_1$ by part a) of Lemma \ref{lem:fundrels}. However, it is unclear a priori that $G_2$ is also equal to $0$. We use the algorithm described in the proof of Proposition \ref{propMCB} to get a sequence of refined basic operations that transforms $G_1$ into $G_2$. We first make the edge of $G_1$ of color $2$ nondisconnecting, and then adjust its position (the third step) and relabel the nodes (last step) to get $G_2$:
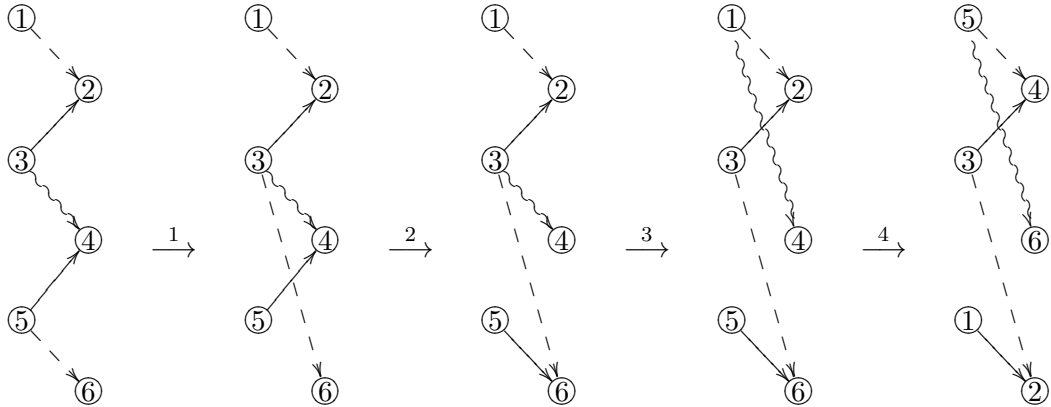
\[
\xymatrix@=15pt{
\ccircle{1} \ar@{-->}[dr] & & & \ccircle{1} \ar@{-->}[dr] & & & \ccircle{1} \ar@{-->}[dr] & & & \ccircle{1} \ar@{-->}[dr] \ar@{~>}[dddr] & & & \ccircle{5} \ar@{-->}[dr] \ar@{~>}[dddr]  & \\
 & \ccircle{2} & & & \ccircle{2}  & &  & \ccircle{2} & &  & \ccircle{2}  & &  & \ccircle{4} \\
\ccircle{3} \ar[ur] \ar@{~>}[dr] & & & \ccircle{3} \ar[ur] \ar@{-->}[dddr] \ar@{~>}[dr] & & & \ccircle{3} \ar[ur] \ar@{~>}[dr] \ar@{-->}[dddr] & & & \ccircle{3} \ar[ur] \ar@{-->}[dddr] & & & \ccircle{3} \ar[ur] \ar@{-->}[dddr]  & \\
 & \ccircle{4} & \overset{1}{\lra} & & \ccircle{4} & \overset{2}{\lra} &  & \ccircle{4} & \overset{3}{\lra} &  & \ccircle{4} & \overset{4}{\lra} &  & \ccircle{6}\\
\ccircle{5} \ar@{-->}[dr] \ar[ur] & & & \ccircle{5} \ar[ur] & &  & \ccircle{5} \ar[dr] & &  & \ccircle{5} \ar[dr] &  &  & \ccircle{1} \ar[dr]& \\
 & \ccircle{6} & & & \ccircle{6} & &  &  \ccircle{6} & &  &  \ccircle{6}  & &  &  \ccircle{2} \\
}
\]
With the notation in the last paragraph of the proof of Proposition \ref{propMCB}, we have $E=(3,4)$ a disconnecting edge, $A_1=\{1,3,5\}$, $B_1=\{2,4,6\}$ a bipartition of the vertex set of $G_1$. We'd like to make $E$ nondisconnecting, with its endpoint in $B_1$ being a leaf. We have
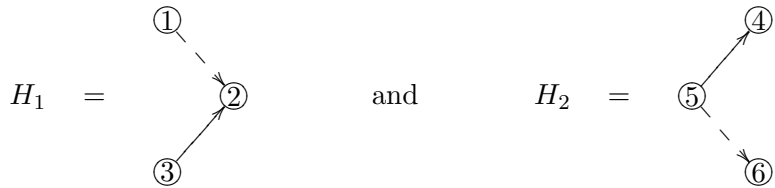
\[\xymatrix@=15pt{
& \ccircle{1} \ar@{-->}[dr] & \\
H_1\quad = & & \ccircle{2} \\
& \ccircle{3} \ar[ur] & \\
}
\quad\quad\quad\quad
\xymatrix@=15pt{
& & \ccircle{4} \\
{\rm{and}\quad\quad\quad\quad H_2\quad =} & \ccircle{5} \ar[ur] \ar@{-->}[dr] & \\
& & \ccircle{6} \\
}
\]
We also have $Y=(5,4)$ of color $y=\xymatrix{\ar[r]&}$ and $O=(5,6)$ of color $o=\xymatrix{\ar@{-->}[r]&}$. The unique vertex $\ccirc{2}$ in $H_1\cap B_1$ is $y$-saturated, and $\ccirc{1}\in H_1\cap A_1$ is $o$-saturated, but $W=\ccirc{3}$ is not $o$-saturated. The refined basic operation involving $W$ and $O$ yields the first transformation.

We now have
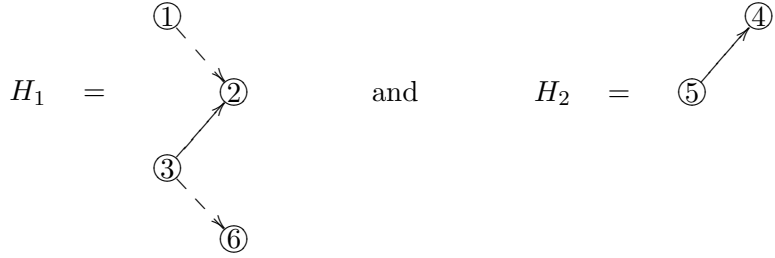
\[\xymatrix@=15pt{
& \ccircle{1} \ar@{-->}[dr] & \\
H_1\quad = & & \ccircle{2} \\
& \ccircle{3} \ar[ur] \ar@{-->}[dr]& \\
& & \ccircle{6} \\
}
\quad\quad\quad\quad
\xymatrix@=15pt{
& & \ccircle{4} \\
{\rm{and}\quad\quad\quad\quad H_2\quad =} & \ccircle{5} \ar[ur] & \\
}
\]
We are in the case $H_2\cap B_1=\{\ccirc{v}\}=\{\ccirc{4}\}$. The edge $E'=(5,4)$ has color $c'=\xymatrix{\ar[r]&}$, different from $c=\xymatrix{\ar@{~>}[r]&}$. $W=\ccirc{6}$ is a vertex in $H_1\cap B_1$ which is not $c'$-saturated, so we can use the refined basic operation involving $E'$ and $W$ as our second transformation, making $E$ a nondisconnecting edge as desired.

We next adjust the position of $E$, in order to get the graph $G_2$. We use the refined operation involving the vertex $\ccirc{1}$ and the edge $(3,4)$. The last transformation involves relabeling the nodes $\ccirc{5}$, $\ccirc{6}$, $\ccirc{2}$, $\ccirc{1}$ and $\ccirc{4}$ by $\ccirc{1}$, $\ccirc{2}$, $\ccirc{4}$, $\ccirc{5}$ and $\ccirc{6}$ respectively.
\end{example}

\begin{corollary}\label{corMCBodd} If $G$ is a canonically oriented MCB-graph of type $(a,a)$, having an odd number of edges, then $G=0$.
\end{corollary}

\begin{proof} Changing the orientation of all the edges of $G$, we obtain a canonically oriented MCB-graph $G'$ of the same type as $G$. It follows from Proposition \ref{propMCB} that $G=G'$. On the other hand, we get by part a) of Lemma \ref{lem:fundrels} that $G'=-G$, hence $G=0$.
\end{proof}

\subsubsection{Step 4} The preceding steps yield the following

\begin{corollary} For $e_{\ll},f_{\ll}$ as in Theorem \ref{thm:main}, the space $(U/F)_{\ll}$ is spanned by MCB-graphs $G_{\mu'}$ of type $\mu'=(a'\geq b')$, with $a'+b'=\min(e_{\ll}+1,r)$ and $b'\geq f_{\ll}$. Moreover, $G_{\mu'}=0$ if $a'=b'$ and $e_{\ll}$ is odd.
\end{corollary}
\begin{proof} The last statement is the content of Corollary \ref{corMCBodd}. We know that $(U/F)_{\ll}$ is spanned by MCB-graphs (Proposition \ref{prop:MCBgenerate}), and the condition $b'\geq f_{\ll}$ follows from the fact that any graph $G$ has at least $\ll_2^i/d_i$ vertices incident to edges of color $i$, and any edge is incident to one vertex in each of the two sets of the bipartition. The number of vertices in the maximal connected component of an MCB-graph of type $\mu'$ is $a'+b'=\min(e_{\ll}+1,r)$.

It remains to show that if $\mu'=(a'\geq b')$, $a'+b'=\min(e_{\ll}+1,r)$ and $b'\geq f_{\ll}$, then there exists an MCB-graph $G_{\mu'}$ of type $\mu'$. Consider $A'$ and $B'$ disjoint sets consisting of $a'$ and $b'$ vertices in $\{\ccirc{1},\cdots,\ccirc{r}\}$ respectively. For every $i=1,\cdots,n$ we draw $\ll^i_2$ edges of color $i$ joining pairs of elements in $A'$ and $B'$, in such a way that no vertex has more than $d_i$ incident edges of color $i$. This is possible since $\ll^i_2/d_i\leq f_{\ll}\leq b'\leq a'$. If the bipartite graph $G$ (with vertex set $A'\cup B'$) obtained in this way is connected, then we get an MCB-graph $G_{\mu'}$ by adding to $G$ the isolated nodes outside $A'\cup B'$. If $G$ is not connected, then it has an edge $E$ of color $c$ contained in a cycle, and a vertex $v$ outside the connected component of $E$. If $v$ is not $c$-saturated, we can move $E$ to make it incident to $v$, and preserve the bipartition of $G$ (as in the refined basic operations of type (1), Remark \ref{remrefbasicop}), thus obtaining a graph with fewer components. If $v$ is $c$-saturated, let $E'$ be an incident edge of color $c$. We move $E$ and $E'$ as in a refined basic operation of type (2), connecting the components of $E$ and $v$. Repeating this procedure will eventually yield a connected graph $G$ and an MCB-graph $G_{\mu'}$ as above.
\end{proof}

\begin{lemma} Consider canonically oriented graphs $G_{\mu'}$ as above, one for each type $\mu'=(a',b')$, with $a'\neq b'$ when $e_{\ll}$ is odd. If
\[\pi=\bigoplus_{\substack{\mu\vdash r\\ \mu=(a\geq b)}}\pi_{\mu}:U\lra\bigoplus_{\substack{\mu\vdash r\\ \mu=(a\geq b)}}U_{\mu}^{\d},\]
then the set $\{\pi(G_{\mu'})\}_{\mu'}$ is linearly independent. In particular, $F=I$ and the graphs $G_{\mu'}$ give a basis of $(U/F)_{\ll}$. This shows that $\rm{dim}((U/I)_{\ll})=m_{\ll}$, where $m_{\ll}$ is as described in Theorem \ref{thm:main}, concluding the proof of our main result.
\end{lemma}
\begin{proof} Note that the number of $G_{\mu'}$'s is precisely $m_{\ll}$, so the last statement follows once we prove the independence of $G_{\mu'}$'s. This is a consequence of the linear independence of $\{\pi(G_{\mu'})\}_{\mu'}$, which in turn follows once we show that for $\mu=(a,b)$, $\mu'=(a',b')$, we have
\begin{enumerate}
\item $\pi_{\mu}(G_{\mu'})=0$ if $b<b'$, and
\item $\pi_{\mu}(G_{\mu'})\neq 0$ if $b=b'$.
\end{enumerate}

Recall that $G_{\mu'}=T_{\mu'}$, for some $n$-tableau $T_{\mu'}$. We have
\[\pi_{\mu}(T_{\mu'})=\sum T_i,\tag{*}\]
where each $T_i$ is an $n$-tableau with entries $1,2$, obtained from a partition $A\sqcup B=\{1,\cdots,r\}$, by setting equal to $1$ and $2$ the entries of $T_{\mu'}$ from $A$ and $B$ respectively.

To prove (1), note that since $|B|=b<b'$, for each $i$ the endpoints of some edge in $G_{\mu'}$ have to be set to the same value, so $T_i$ has repeated entries in some column, i.e. $T_i=0$. It follows that $\pi_{\mu}(G_{\mu'})=\sum T_i=0$.

To prove (2), let $A'\sqcup B'$ be the bipartition of the maximal connected component of $G_{\mu'}$, and take $\mu=(a,b)=(d-b',b')$. The only $n$-tableau(x) $T_i$ in (*) without repeated entries in some column is (are) the $n$-tableau $T_1$ obtained from setting the entries of $A=\{1,\cdots,r\}-B'$ to $1$, and the entries of $B=B'$ to $2$ (and if $|A'|=|B'|$, the $n$-tableau $T_2$ obtained by setting the entries of $A=\{1,\cdots,r\}-A'$ to $1$ and the entries of $B=A'$ to $2$). Since in the latter case $e_{\ll}$ must be even, we get in fact that $T_1=T_2$, since $T_1$ and $T_2$ differ by an even number of transpositions within columns, and by permutations of columns of size $1$. It follows that it's enough to prove that $T_1\neq 0$.

Up to permutations within columns, and permutations of columns of the same size, we may assume that
\[T_1=c_{\ll}\cdot m = c_{\ll}\cdot z_{(A,\cdots,A)}\cdot z_{(B,\cdots,B)},\]
where $A=\{1,\cdots,a\}$ and $B=\{a+1,\cdots,a+b\}$, i.e. $T_1=T_1^1\o\cdots\o T_1^n$, with
\[\Yvcentermath1 T_1^i=\young(11\cdots111\cdots22\cdots,22\cdots2).\]

If $a>b$ and $\s=\t\cdot\t'$, with $\t$ a row permutation and $\t'$ a column permutation of the canonical $n$-tableau $T_{\ll}$ of shape $\ll$, then $\s\cdot m\neq m$, unless $\t'=id$. This shows that the coefficient of $m$ in $T_1$ is a positive number, hence $T_1\neq 0$. If $a=b$, $\s\cdot m=m$ and $\t'\neq id$, then $\t'$ must transpose all the pairs $(1,2)$ in the columns of $T_1$ of size $2$. Since $T_1$ has $e_{\ll}$ (an even number) of such columns, the signature of $\t'$ must be $+1$. It follows again that the coefficient of $m$ in $T_1$ is positive and therefore $T_1\neq 0$.
\end{proof}
\end{proof}

\section*{Acknowledgments} 

I would like to thank David Eisenbud for his guidance, and Tony Geramita, Tony Iarrobino, Joseph Landsberg, Steven Sam, Bernd Sturmfels and Jerzy Weyman for helpful discussions. I would also like to thank Dan Grayson and Mike Stillman for making Macaulay2 (\cite{M2}), which has been used at various stages of this project. This work was partially supported by National Science Foundation Grant No. 0964128.


\end{document}